\documentclass[12pt]{amsart}
\usepackage{amssymb,amsmath,bm}
\def\eps{\varepsilon}
\def\e{{\rm e}}
\def\d{{\rm d}}

\def\I {\mathbb{I}}
\def\R {\mathbb{R}}

\def\N {\mathbb{N}}

\def\H {{\mathfrak A}}
\def\W {{\mathfrak M}}
\def\V {{\mathcal V}}
\def\B {{\mathcal B}}
\def\CC {{\mathfrak C}}
\def\D {{\mathcal D}}
\def\E {{\mathcal E}}
\def\A {{\mathbb A}}
\def\AA {{\mathcal A}}

\def\F {{\mathcal F}}

\def\LL {{\mathcal L}}

\def\M {{\mathcal M}}
\def\PP {{\mathcal P}}
\def\SS {{\mathcal S}}

\def \l {\langle}
\def \r {\rangle}
\def \pt {\partial_t}
\def \ptt {\partial_{tt}}
\def \ps {\partial_s}
\def\ddt {\frac{\d}{\d t}}
\newtheorem{proposition}{Proposition}[section]
\newtheorem{theorem}[proposition]{Theorem}
\newtheorem{corollary}[proposition]{Corollary}
\newtheorem{lemma}[proposition]{Lemma}
\theoremstyle{definition}
\newtheorem{definition}[proposition]{Definition}
\newtheorem{remark}[proposition]{Remark}
\newtheorem{example}[proposition]{Example}
\numberwithin{equation}{section}
\def \au {\rm}
\def \ti {\it}
\def \jou {\rm}
\def \bk {\it}
\def \no#1#2#3 {{\bf #1} (#3), #2.}
\def \eds#1#2#3 {#1, #2, #3.}
\title[A new approach to equations with memory]
{A new approach to equations with memory}

\author[M. Fabrizio, C. Giorgi, V. Pata]
{Mauro Fabrizio, Claudio Giorgi, Vittorino Pata}

\address{Universit\`a di Bologna - Dipartimento di Matematica
\newline\indent
Piazza di Porta San Donato 5, 40127 Bologna, Italy} 
\email{fabrizio@dm.unibo.it}
\address{Universit\`a di Brescia - Dipartimento di Matematica
\newline\indent
Via Valotti 9, 25133 Brescia, Italy} 
\email{claudio.giorgi@ing.unibs.it}
\address{Politecnico di Milano - Dipartimento di Matematica ``F.Brioschi''
\newline\indent
Via Bonardi 9,  20133 Milano, Italy}
\email{vittorino.pata@polimi.it}



\thanks{Work partially supported by the Italian PRIN research project 2006
{\it Problemi a frontiera libera, transizioni di fase e modelli
di isteresi}.}
\begin{document}

%%%%%%%%%%%%%%%%%%%%%%%%%%%%%%%%%%%%%%
\begin{abstract}
We discuss a novel approach to the mathematical analysis
of equations with memory, based on the 
notion of a {\it state}.
This is the initial configuration of the system
at time $t=0$
which can be unambiguously determined by the knowledge
of the dynamics for positive times.
As a model, 
for a nonincreasing convex function $G:\R^+\to\R^+$ such that
$$G(0)=\lim_{s\to 0}G(s)>\lim_{s\to\infty}G(s)>0$$
we consider an abstract version
of the evolution equation 
$$\ptt {\bm u}({\bm x},t)-\Delta\Big[G(0) {\bm u}({\bm x},t) 
+\displaystyle\int_0^\infty G'(s) {\bm u}({\bm x},t-s)\d s\Big]=0$$
arising from linear viscoelasticity.
\end{abstract}
%%%%%%%%%%%%%%%%%%%%%%%%%%%%%%%%%%%%%%

\maketitle

\begin{center}
\begin{minipage}{11cm}
\footnotesize
\tableofcontents
\end{minipage}
\end{center}
\newpage

%%%%%%%%%%%%%%%%%%%%%%%%%%%%%%%%%%%%%%%%%%%%
\section{Preamble}
\label{S1}

\subsection{A general introduction to equations with memory}
Many interesting physical phenomena
(such as viscoelasticity, population dynamics or heat flow in real conductors,
to name some) are modelled
by differential equations which are influenced by
the past values of one or more variables in play: the so-called
{\it equations with memory}. 
The main problem in the analysis of
equations of this kind lies in their nonlocal character,
due to the presence of the memory term (in general, the time
convolution of the unknown function
against a suitable memory kernel).
Loosely speaking,
an evolution equation with
memory has the following formal structure:
\begin{equation}
\label{INTR}
\partial_t w(t)= \F(w(t),w^t(\cdot)),\quad t>0,
\end{equation}
where 
$$w^t(s)=w(t-s),\quad s>0,$$
and $\F$ is some operator
acting on $w(t)$, as well as on the
{\it past values} of $w$ up to the actual time $t$.
The function $w$ is supposed to be known for all $t\leq 0$, where 
it need not solve
the differential equation. Accordingly,
the initial datum has the form
$$
w(t)=w_0(t),\quad t\leq 0,
$$
where $w_0$ is a given function defined on $(-\infty,0]$.

A way to circumvent the intrinsic difficulties posed by the problem
is to (try to) rephrase \eqref{INTR} as an ordinary differential equation 
in some abstract space,
by introducing an auxiliary variable accounting for the past history of $w$,
in order to be in a position to exploit the powerful machinery of the theory
of dynamical systems.
This strategy was devised by C.M.~Dafermos \cite{DAF1}, who, in the context of linear
viscoelasticity, proposed to view 
$w^t$ as an additional variable ruled by its own differential equation,
so translating \eqref{INTR} into a differential system acting on an extended space accounting
for the memory component.

However, when dealing with \eqref{INTR},
what one can actually ``measure" is the function $w(t)$ for $t\geq 0$.
The practical consequences
are of some relevance, since
for a concrete realization
of \eqref{INTR} arising from a specific physical model,
the problem of assigning the initial conditions
is not only of theoretical nature.
In particular, it might happen that two {\it different} initial past histories $w_0$
lead to the {\it same} $w(t)$ for $t\geq 0$. From the viewpoint of the dynamics,
such two different initial past histories
are in a fact indistinguishable. This observation suggests that, rather than the
past history $w^t$, one should employ an alternative variable to describe the initial
state of the system, satisfying
the following natural minimality property: 
$$
\text{\it two different initial states
produce different evolutions $w(t)$ for $t\geq 0$.}
$$ From the philosophical side,
this means that the knowledge of $w(t)$ for all $t\geq 0$ determines in a unique
way the initial state of the problem, the only
object that really influences
the future dynamics.

Of course, the main task is then to determine, if possible,
what is a minimal state associated to \eqref{INTR}.
Unfortunately, a universal strategy is out of reach, and the correct choice 
depends on the particular concrete realization of \eqref{INTR}.
Nonetheless, for a large class of equations with memory, where
the memory contribution enters in the form of a convolution integral
with a nonincreasing positive kernel,
a general scheme seems to be applicable.
In this paper, we discuss an abstract evolution equation
with memory arising from linear viscoelasticity, presenting an approach
which can be easily extended and adapted to many other
differential models
containing memory terms.

\subsection{Plan of the paper}
The goal of the next Section~\ref{S2}, of more physical flavor,
is twofold.
First, we present an overview on materials with hereditary memory, dwelling on
the attempts made through the years to construct mathematical models
accounting for memory effects. Next, we illustrate the physical motivations
leading to the concept of minimal state representation.
In Section~\ref{S3}, we introduce an abstract linear
evolution equation with memory
in convolution form, which, besides its remarkable intrinsic interest,
will serve as a prototype to develop the new approach 
highlighted in the previous sections.
After some notation and preliminary assumptions (Section~\ref{S4}),
we recall the history approach devised by Dafermos (Section~\ref{S5}),
whereas, in Section~\ref{S6}, we make a heuristic
derivation of the state framework,
which will be given a suitable functional formulation  
in the subsequent Section~\ref{S7}
and Section~\ref{S8}. There, we prove the existence of a contraction semigroup,
whose exponential stability is established in Section~\ref{S9}, within
standard assumptions on the memory kernel.
In Section~\ref{S10} and Section~\ref{S11},
we discuss the link between the original equation
and its translated version in the state framework. Finally, in Section~\ref{S12},
we compare the history and the state formulations, showing
the advantages of the latter.
%%%%%%%%%%%%%%%%%%%%%%%%%%%%%%%%%%%%%%%%%%%%

%%%%%%%%%%%%%%%%%%%%%%%%%%%%%%%%%%%%%%%%%%%%
\section{A Physical Introduction}
\label{S2}

\subsection{The legacy of Boltzmann and Volterra}
The problem of the correct modelling
of materials with memory has always
represented a major challenge to mathematicians. The origins of modern  
viscoelasticity and, more generally, of the so-called
hereditary systems traditionally trace back to the works of
Ludwig Boltzmann and Vito Volterra \cite{BOL1,BOL2,VOL1,VOL2}, who  
first introduced the notion of memory in connection with the analysis of  
elastic materials. The key assumption in
the hereditary theory of elasticity can be stated in the following way.

\smallskip
\noindent
{\it For an elastic body
occupying a certain region $\B\subset\R^N$ at rest, the  
deformation of the mechanical system at any point ${\bm x}\in\B$
is a function both of the  
instantaneous stress and of all the past stresses at ${\bm x}$.}

\smallskip
\noindent
In other words, 
calling $\bm u=\bm u({\bm x},t)$ the displacement vector at the point 
${\bm x}\in\B$ at time $t\geq 0$,
the infinitesimal strain tensor
$$\bm \eps=\frac12 
\Big[\nabla \bm u +\nabla \bm u^\top \Big]$$
obeys a constitutive relation of the form
$$
{\bm \eps}({\bm x},t)=\tilde {\bm \eps}({\bm \sigma}({\bm x},t),{\bm \sigma}^t({\bm x},\cdot)),
$$
where
$${\bm \sigma}({\bm x},t)
\quad\text{and}\quad{\bm \sigma}^t({\bm x},s)={\bm \sigma}({\bm x},t-s),$$
with $s>0$, are the stress tensor and its past history at $({\bm x},t)$, respectively.
In the same fashion, the inverse relation can be considered; namely,
$$
{\bm \sigma}({\bm x},t)=
\tilde {\bm \sigma}({\bm \eps}({\bm x},t),{\bm \eps}^t({\bm x},\cdot)).
$$
The above representations allow the appearance of discontinuities at time $t$; 
for instance, ${\bm \eps}({\bm x},t)$ may differ from  
$\lim_{s\to 0}{\bm \eps}^t({\bm x},s)$.
In presence of an external force ${\bm f}={\bm f}({\bm x},t)$,
the related motion equation is given by
$$
\partial_{tt}{\bm u}({\bm x},t)= \nabla\cdot{\bm \sigma}({\bm x},t)+{\bm f}({\bm x},t).
$$

The concept of {\it heredity} was proposed by Boltzmann, essentially in the  
same form later developed by Volterra within a rigorous  
functional setting. However, a more careful analysis  
tells some differences between the two approaches. Quoting \cite{II}, 
``...when speaking of Boltzmann and Volterra, we are facing two  
different scientific conceptions springing from two different  
traditions of classical mathematical physics".

Boltzmann's formulation is focused on hereditary elasticity, 
requiring a fading
{\it initial} strain history 
${\bm \eps}^0({\bm x},\cdot)={\bm \eps}^t({\bm x},\cdot)_{|t=0}$
for every ${\bm x}\in\B$, i.e.,
\begin{equation}
\label{ZEROZERO}
\lim_{s\to\infty} {\bm \eps}^0({\bm x},s)=
\lim_{s\to\infty} {\bm \eps}({\bm x},-s)=0,
\end{equation}
so that, for every fixed $({\bm x},t)$,
$$
{\bm \eps}({\bm x},t)= \int^t_{-\infty} \d {\bm \eps}({\bm x},y),
$$
and assuming the  
linear stress-strain constitutive relation at $({\bm x},t)$ in the  
Riemann-Stieltjes integral form
\begin{equation}
\label{BOL1}
{\bm \sigma}({\bm x},t)= \int^t_{-\infty} {\mathbb G}({\bm x},t-y) 
\d {\bm\eps}({\bm x},y),
\end{equation}
where 
$${\mathbb G}={\mathbb G}({\bm x},s),\quad s>0,$$ 
is a fourth order symmetric tensor (for every fixed $s$), nowadays called
{\it Boltzmann function}. In particular, Boltzmann emphasized
a peculiar behavior of viscoelastic solid materials, named  
{\it relaxation property}:  if the solid is held at a constant strain  
[stress] starting from a given time $t_0\geq0$,
the stress [strain] tends (as $t\to\infty$) to  
a constant value which is ``proportional" to the applied constant   
strain [stress]. Indeed,
if
\begin{equation}
\label{EPST0}
{\bm \eps}({\bm x},t)={\bm \eps}({\bm x},t_0)={\bm \eps}_0({\bm x}),\quad \forall t\geq t_0,
\end{equation}
it follows that
\begin{equation}
\label{Relax}
\lim_{t\to\infty}{\bm \sigma}({\bm x},t)
=\lim_{t\to\infty}\int^{t_0}_{-\infty} {\mathbb G}({\bm x},t-y) 
\d {\bm \eps}({\bm x},y)={\mathbb G}_\infty ({\bm x}){\bm \eps}_0({\bm x}),
\end{equation}
where the {\it relaxation modulus}
$${\mathbb G}_\infty({\bm x})= \lim_{s\to\infty} {\mathbb G}({\bm x},s)$$
is assumed to be positive definite.

Conversely, the general theory devised by Volterra to  
describe the constitutive stress-strain relation is based on the  
Lebesgue representation of linear functionals in the history space. In  
this framework, he stated the fundamental postulates of the elastic  
hereditary action:
\begin{itemize}
\item[$\bullet$]
the principle of invariability of the  
heredity,
\item[$\bullet$] the principle of the closed cycle. 
\end{itemize}
In its simpler 
linear version, the Volterra stress-strain constitutive relation reads
\begin{equation}
\label{VOL1eq}
{\bm \sigma}({\bm x},t)= {\mathbb G}_0({\bm x}){\bm \eps}({\bm x},t)+\int^t_{-\infty}  
{\mathbb G}'({\bm x},t-y){\bm \eps}({\bm x},y)\d y,
\end{equation}
where 
$${\mathbb G}_0({\bm x})=\lim_{s\to 0}{\mathbb G}({\bm x},s)$$
and the {\it relaxation function} ${\mathbb G}'({\bm x},s)$ is the  
derivative with respect to $s$
of the Boltzmann function ${\mathbb G}({\bm x},s)$.
It is apparent that \eqref{VOL1eq}  
can be formally obtained from \eqref{BOL1} by means of an integration by parts, 
provided that \eqref{ZEROZERO} holds true. In which case, the Boltzmann and  
the Volterra constitutive relations are equivalent.
It is also worth noting that if \eqref{EPST0} is satisfied for some $t_0\geq 0$
and
$$
\lim_{t\to\infty}\int^{t_0}_{-\infty} {\mathbb G}'({\bm x},t-y)
{\bm \eps}({\bm x},y)\d y=0,
$$
then, in light of \eqref{VOL1eq},
we recover the stress relaxation property \eqref{Relax}, as
$$
\lim_{t\to\infty}{\bm \sigma}({\bm x},t)= 
{\mathbb G}_0({\bm x}){\bm \eps}_0({\bm x})+\lim_{t\to\infty}\int^t_{t_0}  
{\mathbb G}'({\bm x},t-y){\bm \eps}_0({\bm x})\d y=
{\mathbb G}_\infty ({\bm x}){\bm \eps}_0({\bm x}).
$$

The longterm memory appearing in \eqref{BOL1} and \eqref{VOL1eq} raised
some criticism in the scientific community from
the very beginning, due to the
conceptual difficulty to accept the idea of a past history defined
on an infinite time interval (when even the age of the universe  
is finite!). Aiming to overcome such a philosophical objection, Volterra  
circumvented the problem
in a simple and direct way, assuming  
that the past history vanishes before some time $t_{\rm c}\leq 0$ (say, the creation time). 
Hence, \eqref{VOL1eq} is replaced by
$$
{\bm \sigma}({\bm x},t)= {\mathbb G}_0({\bm x}){\bm \eps}({\bm x},t)+\int_{t_{\rm c}}^t  
{\mathbb G}'({\bm x},t-y){\bm \eps}({\bm x},y)\d y,
$$
and the motion equation becomes the well-known {\it Volterra integro-differential equation}
$$
\partial_{tt}{\bm u}({\bm x},t)= \nabla\cdot{\mathbb G}_0({\bm x})\nabla {\bm  
u}({\bm x},t)+\nabla\cdot\int_{t_{\rm c}}^t {\mathbb G}'({\bm x},t-y)\nabla {\bm  
u}({\bm x},y)\d y +{\bm f}({\bm x},t),
$$
which (besides appropriate boundary conditions)
requires only the knowledge
of the initial data $\bm u({\bm x},t_{\rm c})$ and $\partial_t\bm u({\bm x},t_{\rm c})$.

\begin{remark}
Here, we exploited the equality
$${\mathbb F}{\bm S}={\mathbb F}{\bm S}^\top,$$
which holds for any fourth order symmetric tensor ${\mathbb F}$
and any second order tensor ${\bm S}$.
\end{remark}

\subsection{Further developments: the fading memory principle}
In the thirties, Graffi \cite{GRA1,GRA2} applied Volterra's theory
to electromagnetic materials with memory, successfully  
explaining certain nonlinear wave propagation phenomena occurring in the  
ionosphere.
Nonetheless, the modern theory of materials with memory was developed  
after World War II, when the discovery of new  
materials (e.g., viscoelastic polymers) gave a boost to experimental  
and theoretical researches.
In the sixties, a lot of seminal
papers appeared in the literature, dealing with both linear  
and nonlinear viscoelasticity \cite{BO,CN,DAY1,GRIV,GS,KM,LF}.  
In particular, the thermodynamics of materials with memory provided an  
interesting new field of investigations, mainly because some  
thermodynamic potentials, such as entropy and free energy, are
not unique, even up to an additive constant, and their  
definition heavily depends on the choice of the history space (see,  
for instance, \cite{COL}).

Along the same years, Coleman and Mizel \cite{CM1,CM2} introduced 
a main novelty: the notion of {\it fading memory}.  
Precisely, they considered the Volterra constitutive stress-strain relation  
\eqref{VOL1eq}, with the further assumption that the values of the  
deformation history in the far past produce negligible effects on the  
value of the present stress. In other words, the memory of the material  
is fading in time. Incidentally, this also gave the ultimate answer to
the philosophical question of a memory of  
infinite duration. The fading memory principle is mathematically  
stated by endowing the space ${\mathfrak E}$ of 
initial strain histories ${\bm \eps}^0({\bm x},\cdot)$ with a weighted $L^2$-norm
$$
\|{\bm \eps}^0\|^2_{\mathfrak E}
=\int_\B\int_0^\infty h(s)|{\bm \eps}^0({\bm x},s)|^2\d s\d {\bm x},
$$
where the {\it influence function} $h$ is positive,  
monotone decreasing, and controls the relaxation function ${\mathbb G}'$ 
in the following sense:
$$
\int_\B\int_0^\infty h^{-1}(s)|{\mathbb G}'({\bm x},s)|^2\d s\d {\bm x}<\infty.
$$
The theory of Coleman and Mizel encouraged many  
other relevant contributions in the field,
and was the starting point 
of several improvements in viscoelasticity (see  
\cite{DAY2,FM,GG,GRA4,RHN} and references therein).
On the other hand, as pointed out in \cite{FIC1,FIC2}, the fading  
memory principle turns out to be unable to ensure the well-posedness of the  
full motion equation of linear viscoelasticity
\begin{equation}
\label{full_motion}
\partial_{tt}{\bm u}({\bm x},t)= \nabla\cdot{\mathbb G}_0({\bm x})\nabla {\bm  
u}({\bm x},t)+\nabla\cdot\int_0^\infty {\mathbb G}'({\bm x},s)\nabla {\bm  
u}({\bm x},t-s)\d s +{\bm f}({\bm x},t),
\end{equation}
with known initial data $\bm u({\bm x},0)$, $\partial_t\bm u({\bm x},0)$ and 
$\bm u({\bm x},-s)$, for $s>0$. Moreover,
the arbitrariness of the influence  
function $h$ (for a given ${\mathbb G}'$) reflects into the non uniqueness  
of the history space norm topology.

In order to bypass these difficulties, Fabrizio and coauthors
\cite{DPD1,FGM0,FM0,FM}
moved in two directions. Firstly, they looked for more natural conditions on  
the relaxation function ${\mathbb G}'$, focusing 
on the restriction imposed  by the second law of  
thermodynamics. Secondly, they tried to construct  
an intrinsically defined normed history space. To this end, 
in the spirit
of Graffi's work \cite{GRA3}, they suggested that  
any free energy functional endows the history space with a natural norm  
\cite{FGM1,FGM2}. In this direction, starting from
\cite{BO,DAY1}, many other papers proposed new analytic expressions of  
the maximum and minimum free energies  
\cite{DPD1,DGG,FG1,FG2,GEN}.
The first and well-known expression of the Helmholtz potential in  
linear viscoelasticity is the so called {\it Graffi-Volterra free energy density}
$$
\Psi_{\rm G}({\bm x},t)=\frac12{\mathbb G}_{\infty}({\bm x})|\bm  
\eps({\bm x},t)|^2-\frac12\int_0^\infty \l {\mathbb G}'({\bm x},s) 
[\bm \eps({\bm x},t) 
- \bm \eps^t({\bm x},s)],\bm \eps({\bm x},t) - \bm \eps^t({\bm x},s)\r\d s,
$$
where ${\mathbb G}'$ is negative definite with $s$-derivative  
${\mathbb G}''$ positive semidefinite. With this choice of the energy, 
the asymptotic (exponential) stability of the dynamical  
problem \eqref{full_motion} with ${\bm f}=0$ has been
proved, under the further assumption that, for some $\delta>0$, the 
fourth order symmetric tensor
$${\mathbb G}''({\bm x},s)+\delta {\mathbb G}'({\bm x},s)
$$
is positive semidefinite for (almost) every $s>0$
(see \cite{DAF1,FL,GPR,LZ0,RIV,PAT}). Besides, the
form of $\Psi_{\rm G}$  
suggested the introduction of the new displacement history variable (see  
\cite{DAF1})
$$
\bm \eta^t({\bm x},s)={ \bm u}({\bm x},t) - {\bm u}({\bm x},t-s),\quad s>0,
$$
so that, with reference to \eqref{VOL1eq},
$$
\bm \sigma ({\bm x},t)= {\mathbb G}_\infty({\bm x}){\bm \eps}({\bm x},t)-  
\int_0^\infty{\mathbb G}'({\bm x},s) \nabla\bm \eta^t({\bm x},s)\d s.
$$
Accordingly, the dynamical problem \eqref{full_motion} translates into the system
$$
\begin{cases}
\displaystyle\partial_{tt}{\bm u}({\bm x},t)= \nabla\cdot{\mathbb G}_\infty({\bm x})\nabla {\bm  
u}({\bm x},t)-\nabla\cdot\int_{0}^\infty {\mathbb G}'({\bm x},s)\nabla {\bm  
\eta}^t({\bm x},s)\d s +{\bm f}({\bm x},t),\\
\partial_t\bm \eta^t({\bm x},s)= \partial_t\bm u({\bm x},t)-\partial_s\bm  
\eta^t({\bm x},s),
\end{cases}
$$
which requires the knowledge of the initial data
$\bm u({\bm x},0)$, $\partial_t\bm u({\bm x},0)$ and $\bm \eta^0({\bm x},s)$,
where the initial past history
$\bm \eta^0({\bm x},s)$ is taken in the space ${\mathfrak H}$, dictated by $\Psi_{\rm G}$,
of all functions $\bm \eta=\bm \eta({\bm x},s)$ such that
$$
\|{\bm \eta}\|^2_{\mathfrak H}=-\int_\B\int_0^\infty \l {\mathbb G}'({\bm x},s)  
\nabla{\bm \eta}({\bm x},s),\nabla{\bm \eta}({\bm x},s)\r\d s\d {\bm x}<\infty.
$$

\subsection{The concept of state}
\label{subF}
Unfortunately, a new difficulty arises in connection with the above energetic  
approach. Indeed, depending on the form of  ${\mathbb G}'$,
different (with respect to the almost everywhere
equivalence relation) initial past histories ${\bm \eta}_1,{\bm \eta}_2\in{\mathfrak H}$ could  
produce the same solution to the motion problem \eqref{full_motion}
(clearly, with the same initial data $\bm u({\bm x},0)$
and $\partial_t\bm u({\bm x},0)$). This is  
the case when
\begin{equation}
\label{equiv}
\int_\B\int_0^\infty{\mathbb G}'({\bm x},s+\tau)  
\nabla\big[{\bm \eta}_1({\bm x},s)-{\bm \eta}_2({\bm x},s)\big]\d s\d {\bm x}=0,\quad
\forall\tau>0.
\end{equation}

\begin{remark}
\label{REMNOWAY}
As a consequence,  
there is no way to reconstruct the initial past history $\bm \eta^0({\bm x},s)$
of a given material considered at the initial time $t=0$, 
neither from the knowledge
of the actual state of the system, nor assuming 
to know in advance the future dynamics.
\end{remark}

Noll \cite{NOL} tried to solve the problem by 
collecting all equivalent histories, in the sense of \eqref{equiv}, into  
the same equivalence class, named {\it state} of the material with memory.
Nevertheless, any two different histories in the same  
equivalence class satisfy the relation
$$\|{\bm \eta}_1-{\bm \eta}_2\|_{\mathfrak H}\neq0,$$
which implies that 
the history space ${\mathfrak H}$ is {\it not} a state space (unless each
equivalence class is a singleton) and $\Psi_{\rm G}$ is not a  
state function, as required by thermodynamics.
Further efforts have been made to endow the state space of  
materials with memory with a suitable ``quotient"  
topology, which is typically generated by an  
uncountable family of seminorms \cite{DPD2,GF}. In this direction,  
however, there is no hope to recover a natural norm on the state space.  
Indeed, the main obstacle consists in
handling a space where each element is a set containing an infinite 
number of functions (histories).

A different and more fruitful line of investigations was devised
in \cite{DFG} (see also \cite{GEN}), through the introduction of the notion of
a {\it minimal state}. Drawing the inspiration from
the equivalence relation \eqref{equiv},
the authors called minimal state of the system at time $t$ the function of
the variable $\tau>0$
$$
{\bm\zeta}^t({\bm x},\tau)=-\int_0^\infty{\mathbb G}'({\bm x},\tau+s)\,
\nabla{\bm \eta}^t({\bm x},s)\d s
=-\int_0^\infty{\mathbb G}'({\bm x},\tau+s)
\big[{\bm \eps}({\bm x},t)-{\bm \eps}({\bm x},t-s)\big]\d s.
$$
With this position, the 
stress-strain relation takes the compact form
\begin{equation}
\label{NEWSTRESS}
\bm \sigma ({\bm x},t)= {\mathbb G}_\infty({\bm x})\bm \eps({\bm x},t)+
{\bm\zeta}^t({\bm x},0).
\end{equation}
The difficulties of the previous approaches are
circumvented: the minimal state space is
a function space endowed with a natural weighted 
$L^2$-norm arising from the 
free energy functional 
$$
\Psi_{\rm F}({\bm x},t)=\frac12{\mathbb G}_{\infty}({\bm x})|\bm
\eps({\bm x},t)|^2-\frac12\int_0^\infty \l ({\mathbb G}')^{-1}({\bm x},\tau)
\partial_\tau{\bm\zeta}^t({\bm x},\tau), \partial_\tau{\bm\zeta}^t({\bm x},\tau)\r\d 
\tau,
$$
which involves the minimal state 
representation.

\begin{remark}
As a matter of fact, the history and the state frameworks
are comparable, and the latter is more general (see Section~\ref{S12}
for details).
In particular, as devised in \cite{FL1} 
(see also Lemma~\ref{lemmaSH}), it can be 
shown that 
$$\Psi_{\rm F}({\bm x},t)\leq \Psi_{\rm G}({\bm x},t).$$
\end{remark}

\subsection{The problem of initial conditions}
\label{subPIC}
The classical approach to problems with
memory requires the knowledge of the past history of $\bm u$ at time $t=0$,
playing the role of an initial datum of the problem. This
raises a strong theoretical objection: as
mentioned in Remark~\ref{REMNOWAY}, it is physically
impossible to establish the past history
of $\bm u$ up to time $-\infty$ from measurements of the material
at the actual time, or even assuming the dynamics 
known for all $t>0$. On the other hand, 
in the state formulation one needs to know
the {\it initial} state function 
$${\bm\zeta}^0({\bm x},\tau)={\bm\zeta}^t({\bm x},\tau)_{|t=0},$$
which, as we will see, is
the same as knowing the answer of the stress subject to a constant process
in the time interval $(0,\infty)$; namely, the answer in the future.
At first glance, this appears even more conceptually ambiguous
and technically 
difficult than recovering the past history of $\bm u$. We will show
that it is not so. To this end, let us rewrite 
equation~\eqref{full_motion} in the form
\begin{equation}
\label{MOTION2}
\partial_{tt}{\bm u}({\bm x},t)= \nabla\cdot{\mathbb G}_0({\bm x})\nabla {\bm  
u}({\bm x},t)+\nabla\cdot\int_0^t {\mathbb G}'({\bm x},s)\nabla {\bm  
u}({\bm x},t-s)\d s -{\bm F}_0({\bm x},t)+{\bm f}({\bm x},t),
\end{equation}
having set
\begin{align*}
{\bm F}_0({\bm x},t)&=-\nabla\cdot\int_{0}^\infty {\mathbb G}'({\bm x},t+s)\nabla {\bm  
u}({\bm x},-s)\d s\\
&=\nabla\cdot\big[
{\mathbb G}({\bm x},t)\nabla {\bm u}({\bm x},0)
-{\mathbb G}_\infty({\bm x})\nabla {\bm u}({\bm x},0)-{\bm\zeta}^0({\bm x},t)\big].
\end{align*}
Under, say, Dirichlet boundary conditions, and having a given assignment
of the initial values
$\bm u({\bm x},0)$ and $\partial_t\bm u({\bm x},0)$,
the problem is well-posed, whenever
${\bm F}_0$ is available
(i.e., whenever ${\bm\zeta}^0$ is available). 

\begin{remark}
The function ${\bm F}_0$ is not affected
by the choice of the initial data, nor by the presence of the forcing term
$\bm f$.
Moreover, if the initial past history of $\bm u$ is given, the above relation allows us
to reconstruct ${\bm F}_0$. In this respect, the picture is
at least not worse than before.
\end{remark}

\begin{remark}
It is important to point out that, whereas in the previous approach
the {\it whole} past history of $\bm u$ is required, here, in order
to solve the equation up to any given time $T>0$, the values
of ${\bm F}_0({\bm x},t)$ are needed {\it only} for $t\in[0,T]$.
\end{remark}

Agreed that the hypothesis of an assigned
initial past history of $\bm u$
is inconsistent,
we describe an operative method to construct the function ${\bm F}_0$ 
by means of direct measurements, moving from the observations that materials
with memory (such as polymers) are built by means of specific industrial
procedures. 
Assume that a given material, after its artificial generation
at time $t=0$,
undergoes a process in such a way that
$${\bm u}({\bm x},t)={\bm u}({\bm x},0),\quad\forall t>0.$$
In which case, 
the equality
$$
{\bm\zeta}^t({\bm x},0)={\bm\zeta}^0({\bm x},t)
$$
holds. Indeed,
\begin{align*}
{\bm\zeta}^t({\bm x},0)
&=-\int_0^\infty{\mathbb G}'({\bm x},s)
\big[{\bm \eps}({\bm x},0)-{\bm \eps}({\bm x},t-s)\big]\d s\\
&=-\int_t^\infty{\mathbb G}'({\bm x},s)
\big[{\bm \eps}({\bm x},0)-{\bm \eps}({\bm x},t-s)\big]\d s\\
&=-\int_0^\infty{\mathbb G}'({\bm x},t+s)
\big[{\bm \eps}({\bm x},0)-{\bm \eps}({\bm x},-s)\big]\d s={\bm\zeta}^0({\bm x},t).
\end{align*}
Thus, \eqref{NEWSTRESS} entails the relation
$$
{\bm\zeta}^0({\bm x},t)=
\bm \sigma ({\bm x},t)-{\mathbb G}_\infty({\bm x})\bm \eps({\bm x},0),
$$
meaning that ${\bm\zeta}^0$, and in turn ${\bm F}_0$,
can be obtained by measuring the stress $\bm \sigma ({\bm x},t)$, for all times $t>0$,
of a process frozen at the displacement field ${\bm u}({\bm x},0)$.

Next, we consider the equation of motion \eqref{MOTION2}, relative to a material
generated by the same procedure, but delayed of a time $t_{\rm d}>0$.
For this equation, the corresponding function ${\bm F}_0({\bm x},t)$ is now available.

\begin{remark}
As a matter of fact, ${\bm F}_0({\bm x},t)$ is not
{\it simultaneously} available for all $t>0$. However, since
the first process (which constructs  ${\bm F}_0$) keeps going, at any 
given time $T>0$,
referred to the initial time $t=0$ of the problem under
consideration, we have
the explicit expression of
${\bm F}_0({\bm x},t)$ for all $t\in[0,t_{\rm d}+T]$, which is even more than needed
to solve \eqref{MOTION2} on $[0,T]$.
\end{remark}
%%%%%%%%%%%%%%%%%%%%%%%%%%%%%%%%%%%%%%%%%%%%

%%%%%%%%%%%%%%%%%%%%%%%%%%%%%%%%%%%%%%%%%%%%
\section{An Abstract Equation with Memory}
\label{S3}

\noindent
We now turn to the mathematical aspects of the problem,
developing the state approach for an abstract model equation.

Let $H$ be a separable real Hilbert space, and let $A$
be a selfadjoint strictly positive linear operator on $H$ with compact inverse,
defined on a dense domain
$\D(A)\subset H$.
For $t>0$, we consider the abstract homogeneous linear differential equation
with memory of the second order in time
\begin{equation}
\label{ABS}
\ptt u(t)+A\Big[\alpha u(t)-\int_0^\ell\mu(s)u(t-s)\d s\Big]=0.
\end{equation}
Here, $\alpha>0$, $\ell\in(0,\infty]$
and the {\it memory kernel} 
$$\mu:\Omega=(0,\ell)\to(0,\infty)$$
is a (strictly positive) nonincreasing summable function 
of total mass
$$\int_0^\ell\mu(s)\d s\in(0,\alpha),$$
satisfying the condition (automatically fulfilled if $\ell=\infty$)
$$\lim_{s\to\ell}\mu(s)=0.$$
The dissipativity of the system is
entirely contained in the convolution term,
which accounts for the {\it delay} effects: precisely, finite delay if $\ell<\infty$,
infinite delay if $\ell=\infty$.
The equation is supplemented with the initial conditions given
at the time $t=0$
\begin{equation}
\label{ABS-IC}
\begin{cases}
u(0)=u_0,\\
\pt u(0)=v_0,\\
u(-s)_{|s\in\Omega}=\phi_0(s),
\end{cases}
\end{equation}
where $u_0$, $v_0$ and the function $\phi_0$, 
defined on $\Omega$, are prescribed data.

\begin{remark}
A concrete realization of the abstract equation \eqref{ABS} is obtained
by setting $\Omega=\R^+$, $H=[L^2(\B)]^N$, where $\B\subset\R^N$ is a bounded domain
with sufficiently smooth boundary $\partial\B$, 
and 
$$A=-\Delta\quad\text{with}\quad\D(A)=[H^2(\B)]^N\cap [H^1_0(\B)]^N.$$
In that case, calling
$$\alpha-\int_0^s \mu(\sigma)\d \sigma=G(s)
$$
and
$$u(t)={\bm u}({\bm x},t),\quad {\bm x}\in\B,
$$
the equation reads
$$
\begin{cases}
\displaystyle\ptt {\bm u}({\bm x},t)-\Delta\Big[G(0) {\bm u}({\bm x},t) 
+\int_0^\infty G'(s) {\bm u}({\bm x},t-s)\d s\Big]=0,\\
{\bm u}({\bm x},t)_{|{\bm x}\in\partial\B}=0,
\end{cases}
$$
and rules the evolution of the relative
displacement field ${\bm u}$ in a 
homogeneous isotropic linearly viscoelastic
solid occupying a volume $\B$ at rest
\cite{FM,RHN}.
\end{remark}

Putting $\mu(s)=0$ if $s>\ell$, and defining
\begin{equation}
\label{FZERO}
F_0(t)=\int_0^\ell \mu(t+s)\phi_0(s)\d s,
\end{equation}
equation \eqref{ABS} takes the form
\begin{equation}
\label{FORM}
\ptt u(t)+A\Big[\alpha u(t)-\int_0^t\mu(s)u(t-s)\d s-F_0(t)\Big]=0.
\end{equation}
Introducing the Hilbert space 
$$V=\D(A^{1/2}),$$
with the standard inner product and norm
$$\langle u_1,u_2\rangle_V=\langle A^{1/2}u_1,A^{1/2}u_2\rangle_H,
\quad\|u\|_V=\|A^{1/2}u\|_H,$$
we stipulate the following definition of (weak) solution.

\begin{definition}
\label{EU}
Let $u_0\in V$, $v_0\in H$ and $\phi_0:\Omega\to V$ be such that the corresponding
function $F_0$ given by \eqref{FZERO} fulfills
$$F_0\in L^1_{\rm loc}([0,\infty);V).$$ 
A function
$$u\in C([0,\infty),V)\cap C^1([0,\infty),H)$$
is said to be a solution to the Cauchy
problem \eqref{ABS}-\eqref{ABS-IC}
if 
$$
u(0)=u_0,\quad\pt u(0)=v_0,$$
and the equality
$$\l\ptt u(t),w\r+\alpha\l u(t),w\r_V-\int_0^t\mu(s)\l u(t-s),w\r_V \d s-\l F_0(t),w\r_V=0$$
holds for every $w\in V$
and almost every $t>0$, where $\l\cdot,\cdot\r$ denotes
duality.
\end{definition}
%%%%%%%%%%%%%%%%%%%%%%%%%%%%%%%%%%%%%%%%%%%%

%%%%%%%%%%%%%%%%%%%%%%%%%%%%%%%%%%%%%%%%%%%%
\section{Notation and Assumptions}
\label{S4}

\subsection{Notation}
The symbols $\l\cdot,\cdot\r_X$ and $\|\cdot\|_X$
stand for the inner product and the norm on a generic Hilbert space
$X$, respectively. In particular, for the spaces $H$ and $V$,
we have the well-known norm relations
$$\|u\|_V=\|A^{1/2}u\|_H\geq\sqrt{\lambda_1}\,\|u\|_H,
\quad\forall u\in V,$$
where $\lambda_1>0$ is the first eigenvalue of $A$.
We denote by
$$V^*=\D(A^{-1/2})$$ 
the dual space of $V$,
and by $\l\cdot,\cdot\r$ the duality product between $V^*$ and $V$.
We also recall the equality
$$\|w\|_{V^*}=\|A^{-1/2}w\|_H,\quad\forall w\in V^*.$$
For a nonnegative (measurable) function $\omega$ on $\Omega=(0,\ell)$
and for $p=1,2$, we define
the weighted $L^p$-space of $X$-valued functions
$$L^p_\omega(\Omega;X)=\Big\{\psi:\Omega\to X\,:\,
\int_0^\ell\omega(s)\|\psi(s)\|_X^p\d s<\infty\Big\},
$$
normed by
$$\|\psi\|_{L^p_\omega(\Omega;X)}
=\Big(\int_0^\ell\omega(s)\|\psi(s)\|_X^p\d s\Big)^{1/p}.$$
If $p=2$, this is a Hilbert space endowed with the inner product
$$\l\psi_1,\psi_2\r_{L^2_\omega(\Omega;X)}=
\int_0^\ell\omega(s)\l\psi_1(s),\psi_2(s)\r_X\d s.
$$
Finally, given a generic function $\psi:\Omega\to X$, we denote by $D\psi$ its distributional
derivative.

\smallskip
\noindent
{\it A word of warning.} In order to simplify the notation, if $\psi$ 
is any function on $\Omega$, we agree to interpret $\psi(s)=0$ whenever
$s\not\in\Omega$ (in particular, if $\ell<\infty$, whenever
$s>\ell$).

\subsection{Assumptions on the memory kernel}
As anticipated above, $\mu:\Omega\to(0,\infty)$
is nonincreasing and summable.
Setting
$$M(s)=\int_s^\ell\mu(\sigma)\d\sigma=
\int_0^\ell\mu(s+\sigma)\d\sigma,$$
we require that $M(0)<\alpha$.

\begin{remark}
Note that $M(s)>0$ for every $s\in[0,\ell)$, and $\lim_{s\to\ell}M(s)=0$.
\end{remark}

For simplicity, we will take
\begin{equation}
\label{M0}
\alpha-M(0)=1.
\end{equation}
In addition, we suppose that $\mu$ is absolutely continuous
on every closed interval contained in $\Omega$. In particular, $\mu$
is differentiable almost everywhere in $\Omega$ and $\mu'\leq 0$.
Finally, $\mu$ is assumed to be continuous at $s=\ell$, with $\mu(\ell)=0$, if $\ell<\infty$.
Conversely, if $\ell=\infty$, as $\mu$ is nonincreasing and summable, we automatically
have that $\mu(s)\to 0$ as $s\to\infty$.
In fact, we could consider more general kernels as well, allowing
$\mu$ to have a finite or even a countable number of jumps (cf.\ \cite{CP,PAT}).
However, in this work, we will restrict to the continuous case,
in order not to introduce further technical difficulties.
%%%%%%%%%%%%%%%%%%%%%%%%%%%%%%%%%%%%%%%%%%%%

%%%%%%%%%%%%%%%%%%%%%%%%%%%%%%%%%%%%%%%%%%%%
\section{The History Approach}
\label{S5}

\noindent
An alternative way to look at the equation is to work in the
so-called history space framework, devised by Dafermos
in his pioneering paper \cite{DAF1},
by considering the {\it history} variable
$$\eta^t(s)=u(t)-u(t-s),\quad t\geq 0, s\in\Omega,$$
which, formally, fulfills the problem
$$
\begin{cases}
\pt \eta^t(s)=-\ps\eta^t(s)+\pt u(t),\\
\eta^t(0)=0,\\
\eta^0(s)=u_0-\phi_0(s).
\end{cases}
$$
To set the idea in a precise context,
let us introduce the {\it history space}
$$\M=L^2_\mu(\Omega;V),$$
along with the strongly continuous semigroup $R(t)$ of right translations
on $\M$, namely,
$$
(R(t)\eta)(s)=
\begin{cases}
0 & 0<s\le t,\\
\eta(s-t) & s>t,
\end{cases}
$$
whose infinitesimal generator is the linear operator $T$ defined as 
(cf.\ \cite{GP-Terreni,PAT})
$$T\eta=-D\eta,\quad \D(T)=\{\eta\in{\M}\,:\,D\eta\in\M,\,\eta(0)=0\},$$
where
$\eta(0)=\lim_{s\to 0}\eta(s)$ in $V$.
Then, recalling \eqref{M0}, equation \eqref{ABS} translates
into the differential system in the two variables $u=u(t)$ and $\eta=\eta^t(s)$
\begin{equation}
\label{HIS}
\begin{cases}
\ptt u(t)+ A\Big[u(t) +\displaystyle\int_0^\ell\mu(s)\eta^t(s)\d s\Big]=0,\\
\pt\eta^t=T\eta^t+\pt u(t).
\end{cases}
\end{equation}
Accordingly, the initial conditions \eqref{ABS-IC} turn into
\begin{equation}
\label{HIS-IC}
\begin{cases}
u(0)=u_0,\\
\pt u(0)=v_0,\\
\eta^0=\eta_0,
\end{cases}
\end{equation}
where
\begin{equation}
\label{RELINI}
\eta_0(s)=u_0-\phi_0(s).
\end{equation}
Introducing the {\it extended history space}
$$\W=V\times H\times\M,$$
normed by
$$\|(u,v,\eta)\|_\W^2=\|u\|_V^2+\|v\|_H^2+\|\eta\|_\M^2,$$
problem \eqref{HIS}-\eqref{HIS-IC}
generates a
contraction semigroup $\Sigma(t)$ on $\W$ (see \cite{FL,GP-Terreni,PAT}), 
such that, for every $(u_0,v_0,\eta_0)\in\W$,
$$\Sigma(t)(u_0,v_0,\eta_0)=(u(t),\pt u(t),\eta^t).$$
Moreover, $\eta^t$ has the explicit representation
\begin{equation}
\label{REPETA}
\eta^t(s)=
\begin{cases}
u(t)-u(t-s) & 0<s\le t,\\
\eta_0(s-t)+u(t)-u_0 & s>t.
\end{cases}
\end{equation}
Concerning the relation between \eqref{HIS}-\eqref{HIS-IC}
and the original problem \eqref{ABS}-\eqref{ABS-IC},
the following result holds \cite{GP-Terreni}.

\begin{proposition}
\label{propequiveta}
Let $(u_0,v_0,\eta_0)\in\W$. Then, the first component $u(t)$ of
$\Sigma(t)(u_0,v_0,\eta_0)$ solves
\eqref{ABS}-\eqref{ABS-IC} with
$$F_0(t)=\int_0^\ell \mu(t+s)\big\{u_0-\eta_0(s)\big\}\d s.
$$
\end{proposition}

It is easy to see that $\eta_0\in\M$ implies that
$F_0\in L^\infty(\R^+;V)$.
%%%%%%%%%%%%%%%%%%%%%%%%%%%%%%%%%%%%%%%%%%%%

%%%%%%%%%%%%%%%%%%%%%%%%%%%%%%%%%%%%%%%%%%%%
\section{The State Approach}
\label{S6}

\noindent
An essential drawback of the history approach is that,
for given initial data $u_0$ and $v_0$,
two different initial histories may lead to the same
solution $u(t)$, for $t\geq 0$. Somehow, this is not surprising,
since what really enters in the definition of 
a solution to \eqref{ABS}-\eqref{ABS-IC}, rather than $\phi_0$
(which, by \eqref{RELINI}, is related to the initial history $\eta_0$),
is the function $F_0$, defined in \eqref{FZERO} and appearing
in equation \eqref{FORM}.
Thus, from the dynamical viewpoint, two initial
data $\phi_{01}$ and $\phi_{02}$ should be considered by all means {\it equivalent}
when the corresponding function $F_{01}$ and $F_{02}$ coincide,
due to the impossibility to distinguish
their effects in the future.
On this basis,
it seems natural to devise
a scheme where, rather than $\phi_0$, is the function $F_0$ to appear as
the {\it actual} initial datum accounting for the past history of $u$.

In order to translate this insight into a consistent mathematical theory, it is 
quite helpful to 
see first what happens at a {\it formal level}.
To this aim,
for $t\geq 0$ and $\tau\in\Omega$, we introduce the 
(minimal) {\it state} variable
$$\zeta^t(\tau)=\int_0^\ell\mu(\tau+s)\big\{u(t)-u(t-s)\big\}\d s,$$
which fulfills 
the problem
$$
\begin{cases}
\pt \zeta^t(\tau)=\partial_\tau \zeta^t(\tau)+M(\tau)\pt u(t),\\
\zeta^t(\ell)=0,\\
\zeta^0(\tau)=\zeta_0(\tau),
\end{cases}
$$
having set
$$\zeta_0(\tau)=\int_0^\ell\mu(\tau+s)\{u_0-\phi_0(s)\}\d s
=M(\tau)u_0-F_0(\tau).$$
Accordingly, in light of \eqref{M0}, equation
\eqref{ABS} takes the form
$$
\ptt u(t)+ A\big[u(t) +\zeta^t(0)\big]=0,
$$
where
$$\zeta^t(0)=\lim_{\tau\to 0}\zeta^t(\tau)
=\int_0^\ell\mu(s)\big\{u(t)-u(t-s)\big\}\d s.
$$
Rather than $\zeta^t$, it seems more convenient to consider as a state
the new variable
$$\xi^t(\tau)=-\partial_\tau \zeta^t(\tau)
=-\int_0^\ell\mu'(\tau+s)\big\{u(t)-u(t-s)\big\}\d s,$$
which, in turn, fulfills the problem
$$
\begin{cases}
\pt \xi^t(\tau)=\partial_\tau \xi^t(\tau)+\mu(\tau)\pt u(t),\\
\xi^0(\tau)=\xi_0(\tau),
\end{cases}
$$
where the initial datum $\xi_0$ reads
$$
\xi_0(\tau)
=-\int_0^\ell\mu'(\tau+s)\big\{u_0-\phi_0(s)\big\}\d s
=\mu(\tau) u_0+\int_0^\ell\mu'(\tau+s)\phi_0(s)\d s.
$$
If $\ell<\infty$, we
have also the ``boundary" condition
$$\xi^t(\ell)=0,$$
which comes from the very definition of $\xi^t$.
Since $\zeta^t(\ell)=0$, 
we find the relation
\begin{equation}
\label{FORMEQ}
\int_{\tau_0}^\ell \xi^t(\tau)\d\tau=\zeta^t(\tau_0),\quad\forall \tau_0\in\Omega.
\end{equation}
In particular, in the limit $\tau_0\to 0$,
$$\int_0^\ell \xi^t(\tau)\d\tau=\zeta^t(0).
$$
Therefore, \eqref{ABS}-\eqref{ABS-IC} is
(formally) translated into the system
\begin{equation}
\label{BASE2}
\begin{cases}
\ptt u(t)+ A\Big[u(t) +\displaystyle\int_0^\ell \xi^t(\tau)\d\tau\Big]=0,\\
\pt \xi^t(\tau)=\partial_\tau \xi^t(\tau)+\mu(\tau)\pt u(t),
\end{cases}
\end{equation}
with initial conditions
\begin{equation}
\label{BASE2-IC}
\begin{cases}
u(0)=u_0,\\
\pt u(0)=v_0,\\
\xi^0(\tau)=\xi_0(\tau).
\end{cases}
\end{equation}

\begin{remark}
Observe that the nonlocal character of
\eqref{ABS} is not present in \eqref{BASE2} any longer,
since it is hidden in the new variable $\xi^t$.
\end{remark}

At this point, to complete the project,
two major issues need to be addressed:
\begin{itemize}
\item[$\bullet$] Firstly, we have to
write \eqref{BASE2}-\eqref{BASE2-IC} as a differential equation in
a suitable functional space, providing an existence and uniqueness result.
\smallskip
\item[$\bullet$] Secondly, we have to establish a correspondence (not only formal) 
between the solutions to 
\eqref{BASE2}-\eqref{BASE2-IC} and the solutions to the original
problem \eqref{ABS}-\eqref{ABS-IC}. 
\end{itemize}
%%%%%%%%%%%%%%%%%%%%%%%%%%%%%%%%%%%%%%%%%%%%

%%%%%%%%%%%%%%%%%%%%%%%%%%%%%%%%%%%%%%%%%%%%
\section{The State Space}
\label{S7}

\noindent
The first step to set \eqref{BASE2}-\eqref{BASE2-IC}
in a proper functional framework is to interpret
in a correct way the derivative $\partial_\tau$ appearing in
the second equation of \eqref{BASE2}.
We introduce the new memory kernel 
$$\nu(\tau)=\frac{1}{\mu(\tau)}:\Omega\to[0,\infty),$$
and we put
$$\nu(0)=\lim_{\tau\to 0}\nu(\tau).$$
In view of the assumptions on $\mu$, the function $\nu$ is
continuous and nondecreasing
on $\Omega$, with nonnegative derivative (defined a.e.)
$$\nu'(\tau)=-\frac{\mu'(\tau)}{[\mu(\tau)]^2}.$$
Moreover,
$$\lim_{\tau\to \ell}\nu(\tau)=\infty.$$
Introducing the {\it state space}
$$\V=L^2_\nu(\Omega;V),$$
whose norm is related
to the free energy functional 
$\Psi_{\rm F}$ of Section~\ref{subF},
we consider the strongly continuous semigroup
$L(t)$ of left translations on $\V$, defined by
$$(L(t)\xi)(\tau)=\xi(t+\tau).
$$
It is standard matter to verify
that the infinitesimal generator of $L(t)$ is the linear operator $P$ on $\V$ with
domain
$$\D(P)=\big\{\xi\in\V\,:\,D\xi\in\V,\,\xi(\ell)=0\big\},$$
where $\xi(\ell)=\lim_{\tau\to\ell} \xi(\tau)$ in $V$,
acting as
$$P\xi=D\xi,\quad\forall\xi\in\D(P).$$
Note that, if $\xi\in\D(P)$, then $\|\xi\|_V\in C(\Omega)$.

\begin{remark}
\label{remLINFXI}
If $\ell=\infty$, the condition $\xi(\infty)=0$ is automatically
satisfied whenever $\xi,D\xi\in\V$. Indeed, for every $\tau>1$, the function
$\xi$ is absolutely continuous on $[1,\tau]$ with values in $V$, and
$$\xi(\tau)=\xi(1)+\int_1^\tau D\xi(s)\d s.$$
On the other hand,
\begin{align*}
\int_1^\tau\|D\xi(s)\|_V\d s
&=
\int_1^\tau \sqrt{\mu(s)}\,\sqrt{\nu(s)}\,\|D\xi(s)\|_V\d s\\
&\leq \Big(\int_1^\tau \mu(s)\d s\Big)^{1/2}\Big(\int_1^\tau 
\nu(s)\|D\xi(s)\|_V^2\d s\Big)^{1/2}\\
&\leq \sqrt{M(1)}\,\|D\xi\|_\V.
\end{align*}
Thus, 
$$\lim_{\tau\to\infty}\int_1^\tau D\xi(s)\d s$$
exists in $V$, and so does
$$\xi(\infty)=\lim_{\tau\to\infty}\xi(\tau).$$ 
Since the function 
$$\tau\mapsto\nu(\tau)\|\xi(\tau)\|_V^2$$
is summable
and $\nu(\tau)\to\infty$ as $\tau\to\infty$, it must necessarily be $\xi(\infty)=0$.
Arguing in a similar manner, we see that 
$$\xi,D\xi\in\V\quad\Rightarrow \quad \|\xi\|_V\in L^\infty(\Omega).$$
\end{remark}

\begin{lemma}
\label{lemmaP}
For every $\xi\in\D(P)$,
$$\int_0^\ell\nu'(\tau)\|\xi(\tau)\|^2_V\d\tau<\infty,$$
and
the limit
$$\nu(0)\|\xi(0)\|^2_V=\lim_{\tau\to 0} \nu(\tau)\|\xi(\tau)\|^2_V$$
exists finite (equal to zero if $\nu(0)=0$). Moreover,
\begin{equation}
\label{P}
2\l P\xi,\xi\r_\V=-\nu(0)\|\xi(0)\|^2_V-\int_0^\ell\nu'(\tau)\|\xi(\tau)\|^2_V\d\tau\leq 0.
\end{equation}
\end{lemma}

\begin{proof}
We begin to prove the existence of a sequence $\ell_n\uparrow\ell$ such
that
$$\nu(\ell_n)\|\xi(\ell_n)\|^2_V\to 0.$$
If $\ell=\infty$, this is a direct consequence
of the summability of $\nu\|\xi\|^2_V$. Conversely, if $\ell<\infty$,
for every $\tau<\ell$ we have (recalling that $\xi(\ell)=0$)
$$\nu(\tau)\|\xi(\tau)\|^2_V\leq 
\Big(\int_\tau^{\ell}\sqrt{\nu(s)}\,\|D\xi(s)\|_V\d s\Big)^2
\leq (\ell-\tau)\|D\xi\|_\V^2.
$$
Let now $\eps_n<\ell_n$ be any sequence such that $\eps_n\downarrow 0$.
Then,
\begin{align*}
2\l P\xi,\xi\r_\V
&=\lim_{n\to\infty}\int_{\eps_n}^{\ell_n}\nu(\tau)\frac{\d}{\d\tau}\|\xi(\tau)\|_V^2\d\tau\\
&=\lim_{n\to\infty}\Big\{-
\nu(\eps_n)\|\xi(\eps_n)\|^2_V
-\int_{\eps_n}^{\ell_n}\nu'(\tau)\|\xi(\tau)\|_V^2\d\tau\Big\}.
\end{align*}
Since the limit exists finite, both summands have
the same sign and the integral term is monotone, we conclude that
the limit of the sum equals the sum of the limits, so yielding 
equality \eqref{P}.
We are left to demonstrate the implication
$$\nu(0)=0\quad\Rightarrow \quad \nu(0)\|\xi(0)\|^2_V=0.$$
Indeed,
choosing an arbitrary $\tau_0\in\Omega$, for any $\tau<\tau_0$ we have
\begin{align*}
\nu(\tau)\|\xi(\tau)\|^2_V
&\leq 2\nu(\tau)\|\xi(\tau_0)\|_V^2+2\nu(\tau)\Big(\int_\tau^{\tau_0}
\|D\xi(s)\|_V\d s\Big)^2\\
&\leq 2\nu(\tau)\|\xi(\tau_0)\|_V^2+2\Big(\int_\tau^{\tau_0}
\sqrt{\nu(s)}\,\|D\xi(s)\|_V\d s\Big)^2\\
\noalign{\vskip1mm}
&\leq 2\nu(\tau)\|\xi(\tau_0)\|_V^2 +2(\tau_0-\tau)\|D\xi\|_\V^2.
\end{align*}
Therefore,
$$\limsup_{\tau\to 0}\nu(\tau)\|\xi(\tau)\|^2_V\leq 2\tau_0\|D\xi\|_\V^2,$$
and letting $\tau_0\to 0$, the claim follows.
\end{proof}

The following simple lemma will be needed in the sequel.

\begin{lemma}
\label{xiINT}
Let $\xi\in\V$.
Then, $\xi\in L^1(\Omega;V)$ and
$$\int_0^\ell \|\xi(\tau)\|_V\d\tau\leq \sqrt{M(0)}\,\|\xi\|_\V.
$$ 
As a byproduct, the function
$$t\mapsto \int_t^\ell \xi(\tau)\d\tau$$
belongs to $C([0,\infty),V)$ and vanishes at infinity.
\end{lemma}

\begin{proof}
Using the H\"older inequality,
$$\int_0^\ell \|\xi(\tau)\|_V\d\tau= \int_0^\ell \sqrt{\mu(\tau)}\,
\sqrt{\nu(\tau)}\,\|\xi(\tau)\|_V\d\tau\leq \sqrt{M(0)}\,\|\xi\|_\V,
$$
as claimed.
\end{proof}
%%%%%%%%%%%%%%%%%%%%%%%%%%%%%%%%%%%%%%%%%%%%

%%%%%%%%%%%%%%%%%%%%%%%%%%%%%%%%%%%%%%%%%%%%
\section{The Semigroup in the Extended State Space}
\label{S8}

\noindent
We are now in a position to
formulate \eqref{BASE2}-\eqref{BASE2-IC} as an abstract evolution equation
on a suitable Hilbert space. To this end, we introduce
the {\it extended state space}
$$\H=V\times H\times\V,$$
normed by
$$\|(u,v,\xi)\|_\H^2=\|u\|_V^2+\|v\|_H^2+\|\xi\|_\V^2,$$
and the linear operator $\A$ on $\H$,
with domain
$$\textstyle
\D(\A)=\big\{
(u,v,\xi)\in\H \,:\,
v\in V,\,
u+\int_0^{\ell}\xi(\tau)\d\tau\in \D(A),\,
\xi\in \D(P)\big\},
$$
acting as
$$
\A(u,v,\xi)=
\textstyle\big(
v,
-A\big[u+\int_0^\ell\xi(\tau)\d\tau\big],
P\xi+\mu v\big).
$$
Introducing the 3-component vectors
$$Z(t)=(u(t),v(t),\xi^t)\quad
\text{and}\quad
z=(u_0,v_0,\xi_0)\in\H,$$
we view \eqref{BASE2}-\eqref{BASE2-IC} as the Cauchy problem in $\H$
\begin{equation}
\label{LIN}
\begin{cases}
\displaystyle \ddt Z(t)=\A Z(t),\\
\noalign{\vskip1mm}
Z(0)=z.
\end{cases}
\end{equation}
The following result establishes the existence and uniqueness of a (mild) solution
$$Z\in C([0,\infty),\H).$$

\begin{theorem}
\label{SEMIGROUP}
Problem~\eqref{LIN} generates a 
contraction semigroup $S(t)=\e^{t\A}$ on $\H$
such that
$$Z(t)=S(t)z,\quad \forall t\geq 0.$$
Moreover,
the energy equality
\begin{equation}
\label{ENE}
\ddt\|S(t)z\|^2_\H
=-\nu(0)\|\xi^t(0)\|^2_V-\int_0^\ell\nu'(\tau)\|\xi^t(\tau)\|^2_V\d\tau
\end{equation}
holds for every $z\in\D(\A)$.
\end{theorem}

\begin{proof}
On account of the classical
Lumer-Phillips theorem \cite{PAZ}, we know that $\A$ is the infinitesimal
generator of a contraction semigroup on $\H$ provided that
\begin{itemize}
\item[(i)] the inequality $\l \A z, z\r_{\H}\leq 0$ holds for every $z\in\D(\A)$; and
\smallskip
\item[(ii)] the map $\I-\A:\D(\A)\to\H$ is onto.
\end{itemize}

Concerning point (i), from \eqref{P} we see at once that
$$\l \A z, z\r_{\H}=\l P\xi,\xi\r_\V\leq 0,$$
for every $z=(u,v,\xi)\in\D(\A)$.

In order to prove (ii),
let $z_\star=(u_\star,v_\star,\xi_\star)\in\H$ be given. 
We look for a solution $z=(u,v,\xi)\in\D(\A)$ to the equation
$$(\I-\A)z=z_\star,$$
which, written in components, reads
\begin{equation}
\label{MEDSYS}
\begin{cases}
u-v=u_\star,\\
\displaystyle v+A\Big[u+\int_0^\ell\xi(\tau)\d\tau\Big]=v_\star,\\
\xi(\tau)-D\xi(\tau)-\mu(\tau)v=\xi_\star(\tau).
\end{cases}
\end{equation}
Given a function $g$ on $\Omega$ (extended on the whole real line by setting $g(s)=0$
if $s\not\in\Omega$) and
denoting
$$\E(s)=\e^s\chi_{(-\infty,0]}(s),$$
we consider the convolution product in $\R$ of $\E$ and $g$
at the point $\tau\in\Omega$
$$(\E*g)(\tau)=\int_{\R}\E(\tau-s)g(s)\d s
=\int_\tau^\ell \e^{\tau-s}g(s)\d s.
$$
It is well known (see, e.g., \cite{HS}) that 
$$g\in L^2(\Omega)\quad\Rightarrow\quad \E*g\in L^2(\Omega)$$
and
$$\|\E*g\|_{L^2(\Omega)}\leq \|g\|_{L^2(\Omega)}.$$
For any fixed $v\in V$, we define the function
\begin{equation}
\label{xiGUESS}
\xi(\tau)=v(\E*\mu)(\tau)+(\E*\xi_\star)(\tau).
\end{equation}
We begin to show that $\xi\in\V$.
Indeed,
\begin{align*}
\int_0^\ell \nu(\tau)\|\xi(\tau)\|_V^2\d \tau
&\leq
2\|v\|_V^2\int_0^\ell \nu(\tau)|(\E*\mu)(\tau)|^2\d\tau
+2\int_0^\ell \nu(\tau)\|(\E*\xi_\star)(\tau)\|_V^2\d\tau\\
&\leq
2\|v\|_V^2\int_0^\ell \Big(\int_\tau^\ell\e^{\tau-s}\sqrt{\nu(\tau)}\,
\mu(s)\d s\Big)^2\d\tau\\
&\quad +2\int_0^\ell \Big(\int_\tau^\ell\e^{\tau-s}\sqrt{\nu(\tau)}\,
\|\xi_\star(s)\|_V\d s\Big)^2\d\tau\\
&\leq 2\|v\|_V^2\|\E*\sqrt{\mu}\,\|_{L^2(\Omega)}^2
+2\|\E*(\sqrt{\nu}\,\|\xi_\star\|_V)\|_{L^2(\Omega)}^2\\
\noalign{\vskip1mm}
&\leq 2M(0)\|v\|_V^2
+2\|\xi_\star\|_\V^2.
\end{align*}
Moreover,
$$\|\xi(\tau)\|_V\leq \int_\tau^\ell \e^{\tau-s}
\big\{\mu(s)\|v\|_V+\|\xi_\star(s)\|_V\big\}\d s.
$$
Since (cf.\ Lemma~\ref{xiINT})
$$s\mapsto \mu(s)\|v\|_V+\|\xi_\star(s)\|_V\in L^1(\Omega),$$
we conclude that
$$\lim_{\tau\to\ell}\|\xi(\tau)\|_V=0.$$
Taking the distributional derivative in both sides
of \eqref{xiGUESS}, it is apparent that such a $\xi$ 
satisfies the third equation
of \eqref{MEDSYS}. Moreover, by comparison,
it is readily seen that $D\xi\in\V$.
In summary, $\xi\in\D(P)$ (for any given $v\in V$)
and fulfills the third equation
of \eqref{MEDSYS}.
At this point, we plug $\xi$ into the second equation
of \eqref{MEDSYS}, reading $u$ from the first one.
Noting that
$$\gamma=1+\int_0^\ell\Big(\int_\tau^\ell \e^{\tau-s}\mu(s)\d s\Big)\d\tau
=1+\int_0^\ell \mu(s)\big\{1-\e^{-s}\big\}\d s>0,$$
we obtain
\begin{equation}
\label{LPv}
v+\gamma Av=v_\star-A(u_\star+w),
\end{equation}
having set
$$w=\int_0^\ell (\E*\xi_\star)(\tau)\d\tau.
$$
The elliptic equation \eqref{LPv} admits a (unique) solution $v\in V$,
provided that its right-hand side
belongs to $V^*$, which immediately follows from $w\in V$. Indeed, using the
H\"older inequality,
\begin{align*}
\|w\|_{V}&\leq
\int_0^\ell \|(\E*\xi_\star)(\tau)\|_V\d\tau\\
&\leq
\int_0^\ell (\E*\|\xi_\star\|_V)(\tau)\d\tau\\
&\leq
\int_0^\ell \sqrt{\mu(\tau)}\,\big(\E*(\sqrt{\nu}\,\|\xi_\star\|_V)\big)(\tau)\d\tau\\
&\leq
\sqrt{M(0)}\,\|\xi_\star\|_\V.
\end{align*}
Finally, by comparison, we learn that
$$u+\int_0^\ell\xi(\tau)\d\tau=A^{-1}(v_\star-v)\in\D(A).$$
This completes the proof of point (ii).

We now appeal to a general result of the theory of linear semigroups \cite{PAZ}.
Namely, if $z\in\D(\A)$, then 
$$S(t)z\in\D(\A),\quad \forall t\geq 0,$$
and
$$
\ddt \|S(t)z\|^2_\H
=2\l \A S(t)z, S(t)z\r_{\H}.
$$
On the other hand, since 
$$\l \A S(t)z, S(t)z\r_{\H}=\l P\xi^t,\xi^t\r_\V,$$
the energy equality \eqref{ENE} follows from Lemma~\ref{lemmaP}.
\end{proof}

\begin{corollary}
The third component $\xi^t$ of the solution $S(t)z$ (the state)
has the explicit representation formula
\begin{equation}
\label{REP}
\xi^t(\tau)=\xi_0(t+\tau)+\mu(\tau)u(t)-\mu(t+\tau)u_0+
\int_0^t\mu'(\tau+s) u(t-s)\d s,
\end{equation}
which is valid 
for every $z=(u_0,v_0,\xi_0)\in\H$.
\end{corollary}

\begin{proof}
Assume first that $z$ lies in a more regular space, so that
$\pt u\in L^1_{\rm loc}([0,\infty);V)$. Then, $\xi^t$ satisfies the
nonhomogeneous Cauchy problem in $\V$
$$
\begin{cases}
\displaystyle \ddt \xi^t=P \xi^t+\mu\pt u(t),\\
\noalign{\vskip1mm}
\xi^0=\xi_0.
\end{cases}
$$
Applying the variation
of constants to the semigroup $L(t)=\e^{t P}$ (see \cite{PAZ}),
we obtain
$$
\xi^t(\tau)=\xi_0(t+\tau)
+\int_0^t\mu(t+\tau-s) \pt u(s)\d s.
$$
The desired conclusion \eqref{REP} is drawn integrating by parts.
Using a standard approximation argument, the representation
formula holds for all $z\in\H$.
\end{proof}

\begin{remark}
In the above corollary, the continuity properties of $\mu$ play a crucial role
when integrations by parts occur.
Nonetheless, if $\mu$ has jumps, it is still possible to find 
a representation formula, which contains extra terms accounting for the jumps of $\mu$.
\end{remark}

\begin{remark}
The state variable $\xi^t$ is {\it minimal} in the following sense:
if $(u(t),\pt u(t),\xi^t)$ is a solution to \eqref{LIN} with $u(t)=0$ for every $t\geq 0$,
then $\xi^t$ is identically zero. Indeed, on account of \eqref{LIN} and \eqref{REP},
$$\xi^t(\tau)=\xi_0(t+\tau),\quad\forall t\geq 0,$$
and
$$0=\int_0^\ell\xi^t(\tau)\d\tau=\int_0^\ell\xi_0(t+\tau)\d\tau
=\int_t^\ell\xi_0(\tau)\d\tau,\quad\forall t\geq 0,
$$
which implies that $\xi_0=0$ and, in turn, $\xi^t=0$.
\end{remark}
%%%%%%%%%%%%%%%%%%%%%%%%%%%%%%%%%%%%%%%%%%%%

%%%%%%%%%%%%%%%%%%%%%%%%%%%%%%%%%%%%%%%%%%%%
\section{Exponential Stability}
\label{S9}

\subsection{Statement of the result}
We prove
the exponential stability of the semigroup $S(t)$ on $\H$,
assuming in addition that $\mu$ satisfies
\begin{equation}
\label{SUF}
\mu'(s)+\delta\mu(s)\leq 0,
\end{equation}
for some $\delta>0$ and almost every $s\in\Omega$.

\begin{theorem}
\label{EXPthm}
Let $\mu$ satisfy \eqref{SUF}.
Then, there exist $K>1$ and $\omega>0$ such that
\begin{equation}
\label{EXPformula}
\|S(t)z\|_\H\leq K\|z\|_\H\,\e^{-\omega t},
\end{equation}
for every $z\in\H$.
\end{theorem}

Before proceeding to the proof, some comments are in order.
Condition \eqref{SUF} is quite popular in the literature;
indeed, it has been employed by several authors
to prove the exponential decay of semigroups related
to various equations with memory in the history space framework
(e.g., in connection with the present equation, \cite{FL,GPR,LZ0,RIV}).
On the other hand, the recent paper \cite{PAT} shows
that the exponential decay for such semigroups can be obtained
under the weaker condition
\begin{equation}
\label{MU}
\mu(\sigma+s)\leq C\e^{-\delta \sigma}\mu(s),
\end{equation}
for some $C\geq 1$, every $\sigma\geq 0$ and almost every $s\in\Omega$,
assuming that the set where $\mu'=0$ is not too large (in a suitable
sense).
It is apparent that \eqref{MU} and \eqref{SUF} coincide if $C=1$.
However, if $C>1$, then \eqref{MU} is much more general.
For instance, it is always satisfied
when $\ell<\infty$ (provided that $\mu$ fulfills
the general assumptions of Section~\ref{S4}).
On the contrary, \eqref{SUF} does not allow $\mu$ to have
flat zones, or even horizontal inflection points.
As shown in \cite{CP}, condition \eqref{MU} is actually {\it necessary}
for the exponential decay in the history space framework.
This is true also in the state framework.

\begin{proposition}
Assume that the semigroup $S(t)$ on $\H$ is exponentially stable.
Then, $\mu$ fulfills \eqref{MU}.
\end{proposition}

We omit the proof of the proposition, which can be obtained 
along the lines of \cite{CP},
showing that the exponential stability of $S(t)$ implies the exponential stability
of the left-translation semigroup $L(t)$ on $\V$.

We finally point out that, although we stated the theorem using \eqref{SUF},
the result is still true under the more general hypotheses of \cite{PAT}
(but a much more complicated proof is needed).

\subsection{Proof of Theorem \ref{EXPthm}}
Appealing to the continuity of $S(t)$, it is enough to prove inequality
\eqref{EXPformula} for all $z\in\D(\A)$.
Fix then 
$$z=(u_0,v_0,\xi_0)\in\D(\A),$$
and
denote 
$$S(t)z=(u(t),v(t),\xi^t)\in\D(\A).$$
Introducing the energy
$$E(t)=\frac12\|S(t)z\|_\H^2,$$
and writing \eqref{SUF} in terms of $\nu$ as
$$\nu'(\tau)\geq \delta\nu(\tau),$$
on account of \eqref{ENE}
we derive the differential inequality
\begin{equation}
\label{UNO}
\ddt E(t)\leq-\frac12\int_0^\ell\nu'(\tau)\|\xi^t(\tau)\|^2_V\d\tau
\leq-\frac\delta2\|\xi^t\|_\V^2.
\end{equation}
For an arbitrary $\beta\in\Omega$,
we define the (absolutely continuous) function $\rho:\Omega\to[0,1]$
$$\rho(\tau)=
\begin{cases}
\beta^{-1}\tau& \tau\leq\beta,\\
1&\tau>\beta,
\end{cases}
$$
and we consider the further functionals
\begin{align*}
\Phi_1(t)&=-\int_0^\ell\rho(\tau)\l v(t),\xi^t(\tau)\r_H\d \tau,\\
\Phi_2(t)&=\l v(t),u(t)\r_H.
\end{align*}
Recalling Lemma~\ref{xiINT},
\begin{equation}
\label{BIB}
\int_0^\ell \|\xi^t(\tau)\|_V\d\tau
\leq \sqrt{M(0)}\,\|\xi^t\|_\V.
\end{equation}
Thus, from the continuous embedding $V\subset H$,
\begin{equation}
\label{CTRL}
|\Phi_\imath(t)|\leq c_0 E(t),\quad \imath=1,2,
\end{equation}
for some $c_0>0$ independent of the choice of $z\in\D(\A)$.

\begin{lemma}
\label{lemmaPHI1}
There is $c_1>0$ independent of $z$ such that
\begin{equation}
\label{DUE}
\ddt\Phi_1(t)\leq M(\beta)\Big\{\frac1{12}\|u(t)\|_V^2-\frac12 \|v(t)\|_H^2
+c_1\|\xi^t\|_\V^2\Big\}.
\end{equation}
\end{lemma}

\begin{proof}
We have
$$\ddt\Phi_1=-\int_0^\ell\rho(\tau)\l \pt v,\xi(\tau)\r_H\d \tau
-\int_0^\ell\rho(\tau)\l v,\pt \xi(\tau)\r_H\d \tau.$$
We now estimate the two terms of the right-hand side, 
exploiting the equations of \eqref{LIN} and 
the integral control \eqref{BIB}.
For the first one,
\begin{align*}
-\int_0^\ell\rho(\tau)\l \pt v,\xi(\tau)\r_H\d \tau
&=\int_0^\ell\rho(\tau)\l u,\xi(\tau)\r_V\d \tau
+\int_0^\ell\rho(\tau)
\Big(\int_0^\ell\l\xi(\tau'),\xi(\tau)\r_V\d\tau'\Big)\d \tau\\
&\leq\|u\|_V\int_0^\ell\|\xi(\tau)\|_V\d \tau
+\Big(\int_0^\ell\|\xi(\tau)\|_V\d \tau\Big)^2\\
&\leq\sqrt{M(0)}\,\|u\|_V\|\xi\|_\V+M(0)\|\xi\|_\V^2\\
&\leq\frac1{12} M(\beta)\|u\|_V^2
+M(0)\Big(1+\frac{3}{M(\beta)}\Big)\|\xi\|_\V^2.
\end{align*}
Concerning the second term, 
we preliminarily observe that,
since $\xi\in\D(P)$ for all times, we have (cf.\ Remark~\ref{remLINFXI})
$$\sup_{\tau\in\Omega}\|\xi(\tau)\|_V<\infty\quad\text{and}\quad
\|\xi(\ell)\|_V=0.$$
Thus, an integration by parts gives
$$-\int_0^\ell\rho(\tau)\l v,P\xi(\tau)\r_H\d \tau
=-\int_0^\ell\rho(\tau)\frac{\d}{\d\tau}\l v,\xi(\tau)\r_H\d \tau
=\frac1\beta\int_0^\beta\l v,\xi(\tau)\r_H\d \tau.
$$
Hence,
\begin{align*}
-\int_0^\ell\rho(\tau)\l v,\pt \xi(\tau)\r_H\d \tau
&=-\Big(\int_0^\ell\rho(\tau)\mu(\tau)\d\tau\Big)\|v\|_H^2
+\frac1\beta\int_0^\beta\l v,\xi(\tau)\r_H\d \tau\\
&\leq-M(\beta)\|v\|_H^2
+\frac{1}{\beta\sqrt{\lambda_1}\,}\|v\|_H\int_0^\ell\|\xi(\tau)\|_V\d\tau\\
&\leq-M(\beta)\|v\|_H^2
+\frac{\sqrt{M(0)}\,}{\beta\sqrt{\lambda_1}\,}\|v\|_H\|\xi\|_\V\\
&\leq-\frac12 M(\beta)\|v\|_H^2
+\frac{M(0)}{2\beta^2\lambda_1 M(\beta)}\|\xi\|_\V^2.
\end{align*}
Collecting the above inequalities, the conclusion follows.
\end{proof}

\begin{lemma}
\label{lemmaPHI2}
The functional $\Phi_2(t)$ fulfills the differential inequality
\begin{equation}
\label{TRE}
\ddt\Phi_2(t)\leq -\frac34\|u(t)\|_V^2+\|v(t)\|^2_H+M(0)\|\xi^t\|_\V^2.
\end{equation}
\end{lemma}

\begin{proof} 
By virtue of \eqref{LIN} and \eqref{BIB},
\begin{align*}
\ddt\Phi_2
&=-\|u\|^2_V+\|v\|^2_H-\int_0^\ell\l u,\xi(\tau)\r_V\d \tau\\
&\leq-\|u\|^2_V+\|v\|^2_H
+\sqrt{M(0)}\,\|u\|_V\|\xi\|_\V\\
&\leq-\frac34\|u\|^2_V+\|v\|^2_H
+M(0)\|\xi\|_\V^2,
\end{align*}
as claimed.
\end{proof}

At this point, we define the functional
$$\Phi(t)=\frac{3}{M(\beta)}\Phi_1(t)+\Phi_2(t),
$$
which, due \eqref{DUE} and \eqref{TRE}, satisfies the differential
inequality
\begin{equation}
\label{QUATTRO}
\ddt\Phi(t)+E(t)\leq c_2\|\xi^t\|_\V^2,
\end{equation}
for some $c_2>0$ independent of $z$.
Besides, in light of \eqref{CTRL},
\begin{equation}
\label{CTRLPHI}
|\Phi(t)|\leq c_3 E(t),
\end{equation}
with $c_3=c_0(3/M(\beta)+1)$.
Finally,
we fix
$$\eps=\min\Big\{\frac{\delta}{2c_2},\frac{1}{2c_3}\Big\}
$$
and we set
$$\Psi(t)=E(t)+\eps\Phi(t).$$
Note that, by \eqref{CTRLPHI},
$$
\frac12E(t)\leq\Psi(t)\leq\frac32 E(t),
$$
and in turn, by \eqref{UNO} and \eqref{QUATTRO},
$$
\ddt \Psi(t)+2\omega\Psi(t)\leq 0,
$$
with $\omega=\eps/3$.
Therefore, the standard Gronwall lemma yields
$$\|S(t)z\|_\H^2=2E(t)\leq4\Psi(t)\leq 4\Psi(0)\e^{-2\omega t}
\leq 6E(0)\e^{-2\omega t}= 3\|z\|_\H^2\,\e^{-2\omega t}.
$$
The proof of Theorem~\ref{EXPthm} is completed.

\begin{remark}
Observe that the proof of Theorem \ref{EXPthm} is carried out 
employing {\it only}
energy functionals, and it makes no use of
linear semigroup techniques. Thus,
the same energy functionals
can be exploited to analyze semilinear versions of the problem
(for instance, to prove the existence of absorbing sets and global attractors). 
\end{remark}
%%%%%%%%%%%%%%%%%%%%%%%%%%%%%%%%%%%%%%%%%%%%

%%%%%%%%%%%%%%%%%%%%%%%%%%%%%%%%%%%%%%%%%%%%
\section{The Original Equation Revisited}
\label{S10}

\noindent
Somehow, this novel state approach 
urges us to consider the original problem under a different
perspective.
Indeed, as we saw in Section~\ref{S6},
the solutions to \eqref{ABS}-\eqref{ABS-IC}
are determined, besides by $u_0$ and $v_0$, by the knowledge of the 
function $F_0$,
and not by the particular form of the initial past history $\phi_0$.
Therefore, with reference to Definition~\ref{EU},
we introduce the class of {\it admissible past history functions} 
$$\AA=\Big\{\phi:\Omega\to V\,:\,
t\mapsto\int_0^\ell \mu(t+s)\phi(s)\d s\in L^1_{\rm loc}([0,\infty);V)\Big\},
$$
and we define the linear map 
$$\Lambda:\AA\to L^1_{\rm loc}([0,\infty);V)$$
as
$$\phi\mapsto \Lambda \phi(t)=\int_0^\ell \mu(t+s)\phi(s)\d s.
$$
Note that $\Lambda\phi(t)=0$ if $t\geq\ell$.
Accordingly, we define the class of {\it state functions}
$$\SS=\Lambda\AA.$$
Clearly (and this is really the point), the map $\Lambda$ may not be
injective, meaning that different $\phi\in\AA$ may
lead to the same element of $\SS$.

Coming back to Definition~\ref{EU}, the assumption on $F_0$
can now be rephrased as
$$F_0=\Lambda\phi_0\quad\text{with}\quad\phi_0\in\AA,$$
and we can reformulate the definition of solution
to \eqref{ABS} in the following more convenient (and certainly more physical) way.

\begin{definition}
\label{EUNEW}
Let the triplet
$$(u_0,v_0,F_0)\in V\times H\times\SS$$
be given.
A function 
$$u\in C([0,\infty),V)\cap C^1([0,\infty),H)$$
is said to be a solution to equation
\eqref{ABS} with {\it initial state} $(u_0,v_0,F_0)$ if
$$u(0)=u_0,\quad\pt u(0)=v_0,$$
and the equality
$$\l\ptt u(t),w\r+\alpha\l u(t),w\r_V-\int_0^t\mu(s)\l u(t-s),w\r_V \d s-\l F_0(t),w\r_V=0$$
holds for every $w\in V$
and almost every $t>0$.
\end{definition}

In this definition, the initial datum $\phi_0$ has
completely disappeared, since the state function $F_0$ contains all the
necessary information on the past history of the variable $u$
needed to capture the future dynamics of the equation.
Hence, we removed the (unphysical) ambiguity caused
by two different initial histories leading to the same state function,
which, as we saw, is what really enters in the definition of a solution.

\begin{remark}
We point out that the function $F_0(t)$ is not
influenced by the dynamics for $t\geq 0$, nor by the presence of
a possible external force. As a matter of fact, if the initial past history
$\phi_0$ is known, then $F_0$ is uniquely determined by \eqref{FZERO}.
On the other hand,
even if the particular $\phi_0$ leading to $F_0$ is
unknown, in principle, $F_0$ can still be determined
(cf.\ Section~\ref{subPIC}).
\end{remark}

The remaining of the section is devoted to investigate the properties
of the space $\SS$. We begin with a lemma, which provides a precise formulation
of the formal equality \eqref{FORMEQ}, devised in Section~\ref{S6}.

\begin{lemma}
\label{CRUCIAL}
Whenever $\phi\in\AA$, the map
$$\tau\mapsto \int_0^\ell\mu'(\tau+s)\|\phi(s)\|_V\d s$$
belongs to $L^1(t,\infty)$ for every $t>0$, and the equality
\begin{equation}
\label{TEMPO}
\Lambda\phi(t)=-\int_t^\ell \Big(\int_0^\ell\mu'(\tau+s)\phi(s)\d s\Big)\d\tau
\end{equation}
holds for every $t>0$. Moreover, if
$\phi\in L^1_\mu(\Omega;V)$,
then $\phi\in\AA$ and \eqref{TEMPO} holds for every $t\geq 0$.
\end{lemma}

\begin{proof}
Let $\phi\in\AA$ be given. For every fixed $t>0$,
$$\Lambda\phi(t_0)\in V,\quad\text{for some $t_0\leq t$}.$$
Since $\mu$ is a nonincreasing function
and $\Lambda\phi(t_0)$ is a Bochner integral,
this is the same as saying that
$$\int_0^\ell \mu(t+s)\|\phi(s)\|_V\d s
\leq \int_0^\ell \mu(t_0+s)\|\phi(s)\|_V\d s<\infty.$$
Exploiting the equality
$$\mu(t+s)=-\int_t^\ell \mu'(\tau+s)\d\tau,$$
and exchanging the order of integration,
we conclude that
$$\int_0^\ell \mu(t+s)\|\phi(s)\|_V\d s
=-\int_t^\ell \Big(\int_0^\ell\mu'(\tau+s)\|\phi(s)\|_V\d s\Big)\d\tau<\infty.
$$
Hence, 
$$\tau\mapsto \int_0^\ell\mu'(\tau+s)\|\phi(s)\|_V\d s\in L^1(t,\infty),$$
and \eqref{TEMPO} follows from the Fubini theorem.
Concerning the last assertion, just note that $\phi\in L^1_\mu(\Omega;V)$ if and only if
$\Lambda\phi(0)\in V$.
\end{proof}

\begin{remark}
\label{rembel}
As a straightforward consequence of the lemma,
$$\SS\subset C_0([t,\infty),V),\quad\forall t>0,$$
where $C_0$ denotes the space of continuous functions
vanishing at infinity.
\end{remark}

Given $F\in\SS$,
it is then interesting to see what happens to $F(t)$ in the limit
$t\to 0$.
Three mutually disjoint situations may occur:
\smallskip
\begin{itemize}
\item[(i)] $\lim_{t\to 0} F(t)$ exists in $V$;
\smallskip
\item[(ii)] $F\in L^\infty(\R^+;V)$
but $\lim_{t\to 0} F(t)$ does not exist in $V$;
\smallskip
\item[(iii)] $\|F(t)\|_V$ is unbounded in a neighborhood of $t=0$.
\end{itemize}
\smallskip
As we will see, (i) is the most
interesting case in view of our scopes. For this reason, we 
introduce the further space
$$\SS_0=\big\{F\in\SS\,:\,
\exists\,\lim_{t\to 0}F(t)\text{ in } V\big\}.
$$
In light of Remark~\ref{rembel}, it is apparent that
$$\SS_0\subset C_0([0,\infty),V).$$
We preliminary observe that if $F=\Lambda\phi$
with $\phi\in L^1_\mu(\Omega;V)$,
then Lemma~\ref{CRUCIAL}
yields at once
$$\lim_{t\to 0}F(t)=\Lambda\phi(0)\quad\text{in }V,$$
so that $F\in\SS_0$.
However, the picture can be more complicated.
Indeed, it may happen that $F\in\SS_0$ but $\Lambda\phi(0)$
is not defined for {\it any}
$\phi\in \Lambda^{-1}F$, as the following example shows.

\begin{example}
\label{EXINJ}
Consider the kernel
$$\mu(s)=1-s,\quad\Omega=(0,1).
$$
Given any nonzero vector $u\in V$, set
$$F(t)=
\Big[(1-t)\sin 1-\int_t^1\sin\frac1x\d x\Big]\chi_{[0,1]}(t)\,u,
$$
which clearly satisfies
$$\lim_{t\to 0}F(t)=\Big[\sin 1-\int_0^1\sin\frac1x\d x\Big]u.
$$
Then, $F=\Lambda\phi$ with
$$\phi(s)=-\Big[\frac{1}{(1-s)^2}\cos \frac{1}{1-s}\Big] u,$$
but $\Lambda\phi(0)$ is not defined, since $\phi\not\in L^1_\mu(\Omega;V)$.
To complete the argument, we show that, for this particular
kernel, the linear map $\Lambda$ is injective.
Indeed, let $\tilde\phi\in\AA$ be such that
$$0=\Lambda\tilde\phi(t)=
\int_0^{1-t}(1-t-s)\tilde\phi(s)\d s.
$$
The above equality readily implies that $\tilde\phi=0$.
\end{example}

Let us provide examples also for (ii) and (iii).
Again, $u\in V$ is any nonzero vector. 

\begin{example}
With $\mu$ as in the previous example, set
$$F(t)=
\Big[\sin\frac1t-\sin1+t\cos 1-\cos 1\Big]\chi_{[0,1]}(t)\,u.
$$
Then, $F=\Lambda\phi$ with
$$\phi(s)=\Big[\frac{2}{(1-s)^3}\cos \frac{1}{1-s}
-\frac{1}{(1-s)^4}\sin \frac{1}{1-s}\Big] u.$$
Note that $\|F\|_V\in L^\infty(\R^+)$
but $\lim_{t\to 0} F(t)$ does not exist in $V$.
\end{example}

\begin{example}
\label{L1STRICT}
Consider the kernel
$$\mu(s)=\sqrt{\frac{1-s}{s}},\quad\Omega=(0,1),
$$
and set
$$F(t)=
\Big[\int_0^{1-t}\sqrt{\frac{1-t-s}{s(t+s)}}\,\,\d s\Big]\chi_{[0,1]}(t)\,u.
$$
Then, $F=\Lambda\phi$ with
$$\phi(s)=\frac{1}{\sqrt{s}}\,u.$$
It is easily verified that $\|F\|_V$ is summable on $\R^+$
but
$$\lim_{t\to 0}\|F(t)\|_V=\infty.$$
\end{example}

\begin{remark}
A related (and very challenging) question is the following:
given a function 
$F\in C_0([t,\infty),V)$ for every $t>0$,
find (easy to handle)
conditions ensuring that $F\in\SS$.
In fact, the answer seems to be strongly dependent on the
particular choice of the kernel.
For instance, with $\mu$ as in Example~\ref{EXINJ},
$F\in\SS$ if and only if
$F(t)=DF(t)=0$ for $t\geq 1$,
$F\in L^1(\R^+;V)$, and $DF$ is absolutely continuous from
$[t,1]$ into $V$ for all $t>0$.
In which case, $F=\Lambda\phi$ with $\phi(s)=D^2F(1-s)$.
On the other hand, if $\mu(s)=\e^{-s}$, then
$F\in\SS$ if and only if $F(t)=\e^{-t}u$ with $u\in V$
(cf.\ the next Example~\ref{EEXP}).
\end{remark}
%%%%%%%%%%%%%%%%%%%%%%%%%%%%%%%%%%%%%%%%%%%%

%%%%%%%%%%%%%%%%%%%%%%%%%%%%%%%%%%%%%%%%%%%%
\section{Proper States: Recovering the Original Equation}
\label{S11}

\noindent
The purpose of this section is to establish the link between \eqref{LIN}
and the original equation \eqref{ABS}, up to now only formal.
To this end, we have to recall the particular
form of the initial datum $\xi_0$, obtained in a somewhat heuristic way
in Section~\ref{S6}.
This gives a clue that not all the states are apt
to describe the behavior of the original equation, but only
certain particular states having a well defined structure.

\begin{definition}
A vector $\xi\in\V$ is said to be a {\it proper state} if
$$\xi(\tau)=DF(\tau),$$
for some $F\in\SS$.
We denote by 
$\PP$
the normed subspace of $\V$ (with the norm inherited by $\V$)
of proper states.
\end{definition}

For any given kernel $\mu$,
an immediate example of proper state is
$$\xi(\tau)=\mu(\tau)u,\quad u\in V.$$
Indeed, $\xi=DF$ with 
$$F(t)=-M(t)u=-\int_0^\ell \mu(t+s)u\,\d s.$$

\begin{lemma}
\label{LEMMAPROPER}
Let $\xi\in \PP$. Then, there exists a unique $F\in\SS$ such that $\xi=DF$.
Besides, $F$ belongs to $\SS_0$. Moreover, for every $\phi\in\AA$
such that $F=\Lambda\phi$, it follows that
$$\xi(\tau)=\int_0^\ell\mu'(\tau+s)\phi(s)\d s.$$
Conversely, if $\xi\in\V$ has the above representation for some $\phi\in\AA$,
then $\xi\in\PP$ and
$$\xi(\tau)=D\Lambda\phi(\tau).$$
\end{lemma}

\begin{proof}
Let $F\in\SS$ be such that $\xi=DF$. From Lemma~\ref{xiINT},
$DF\in L^1(\Omega;V)$. Hence, the map
$$t\mapsto -\int_t^\ell DF(\tau)\d\tau=F(t)$$
belongs to $C_0([0,\infty),V)$. Therefore, $F\in\SS_0$, and
it is apparent that $F$ is uniquely determined by $DF$.
The remaining assertions
follow by \eqref{TEMPO}.
\end{proof}

Here is a concrete application of the lemma.

\begin{example}
Let $\mu$ and $u$ as in Example~\ref{L1STRICT}, and define
$$F(t)=\Big[\int_0^{1-t}\sqrt{\frac{1-t-s}{s(t+s)}}\,\,
\sin \frac{1}{s}\,\d s\Big]\chi_{[0,1]}(t)\,u.
$$
Then, $F=\Lambda\phi$ with
$$\phi(s)=\Big[\frac{1}{\sqrt{s}}\sin \frac{1}{s}\Big] u.$$
Setting 
$$\xi(\tau)=\int_0^\ell\mu'(\tau+s)\phi(s)\d s
=-\Big[\frac12\int_0^{1-\tau}\frac{1}{\sqrt{s(1-\tau-s)}\,(\tau+s)^{3/2}}\,
\sin \frac{1}{s}\,\d s\Big]u,$$
it is not hard to verify that $\xi\in\V$. From Lemma~\ref{LEMMAPROPER},
we conclude that $\xi=DF\in\PP$.
Note that, as expected, $F\in\SS_0$. Indeed,
$$\lim_{t\to 0}F(t)=\varkappa u\quad\text{in }V,$$
with
$$\varkappa=\lim_{N\to\infty}\int_1^N
\frac1x\sqrt{\frac{x-1}{x}}\,\sin x\,\d x\sim 0.28.
$$
\end{example}

In particular, Lemma~\ref{LEMMAPROPER} says that the map 
$$\Gamma:\PP\to\SS$$
defined as
$$\Gamma\xi(t)=-\int_t^\ell\xi(\tau)\d\tau
$$
is injective. Since 
$$\Gamma\PP\subset\SS_0,$$
and the inclusion $\SS_0\subset\SS$ can be strict,
the map $\Gamma$ is not, in general, onto.
In fact, the inclusion $\Gamma\PP\subset\SS_0$ can be strict either.

\begin{example}
Let
$$\mu(s)=1-s,\quad\Omega=(0,1).$$
Given any nonzero vector $u\in V$, consider the
function
$$F(t)=\big(\sqrt{t}\,-1\big)^2\chi_{[0,1]}(t)\,u.$$
Then, $F=\Lambda\phi$ with
$$\phi(s)=\Big[\frac{1}{2(1-s)^{3/2}}\Big]u.$$
We conclude that $F\in\SS_0$.
On the other hand,
$$DF(\tau)=
\Big[\frac{\sqrt{\tau}\,-1}{\sqrt{\tau}}\Big]u,
$$
which does not belong to $\V$.
\end{example}

We have now all the ingredients to state the main result of the
section.

\begin{theorem}
\label{profCFR}
Let $(u_0,v_0,F_0)\in V\times H\times\SS$. Assume in addition that
$$F_0\in\Gamma\PP.$$
Then, a function $u$ is a solution to
\eqref{ABS} with initial state $(u_0,v_0,F_0)$ (according to Definition~\ref{EUNEW})
if and only if
$$(u(t),\pt u(t),\xi^t)=S(t)(u_0,v_0,\xi_0),$$
with $\xi^t$ as in \eqref{REP} with
$$\xi_0(\tau)
=\mu(\tau) u_0+DF_0(\tau).$$
Conversely, if $u$ is
a solution to
\eqref{ABS} with initial state $(u_0,v_0,F_0)$
and $F_0\not\in\Gamma\PP$, then there is
no corresponding solution in the extended state space.
\end{theorem}

\begin{proof}
Since
$u\in C([0,\infty),V)$, arguing as in the proof of Lemma~\ref{CRUCIAL},
the equality
$$\int_0^\ell\Big(\int_0^t\mu'(\tau+s)u(t-s)\d s\Big)\d\tau=-\int_0^t\mu(s)u(t-s)\d s$$
holds for every $t> 0$.
Thus, using \eqref{REP}, keeping in mind the particular form of $\xi_0$
and the fact that $F_0\in\SS_0$, we readily 
get
\begin{equation}
\label{PPPPPPP}
\int_0^\ell\xi^t(\tau)\d\tau
=M(0)u(t)-\int_0^t\mu(s)u(t-s)\d s
-F_0(t).
\end{equation}
This equality, in light of \eqref{FORM}, \eqref{M0} and \eqref{LIN}, 
proves the first statement.
 
To prove the converse, assume
that $u(t)$ is at the same time a solution to
\eqref{ABS} with initial state $(u_0,v_0,F_0)$,
and equal to the first
component of $S(t)(u_0,v_0,\xi_0)$, for some $\xi_0\in\V$.
We reach the conclusion by showing that $F_0\in\Gamma\PP$.
Indeed, calling now $\xi^t$ the third component
of $S(t)(u_0,v_0,\xi_0)$, from \eqref{FORM}, \eqref{M0} and \eqref{LIN}, 
we obtain again \eqref{PPPPPPP}. Since, by \eqref{REP},
$$
\int_0^\ell\xi^t(\tau)\d\tau
=\int_t^\ell\xi_0(\tau)\d\tau+M(0)u(t)-M(t)u_0-\int_0^t\mu(s)u(t-s)\d s,
$$
we conclude that
$$
\int_t^\ell\big[\mu(\tau)u_0-\xi_0(\tau)\big]\d\tau
=M(t)u_0-\int_t^\ell\xi_0(\tau)\d\tau
=F_0(t).
$$
Hence, 
$$-\mu(\tau) u_0+\xi_0(\tau)=DF_0(\tau),$$
meaning that $\xi_0-\mu u_0\in\PP$ and $F_0=\Gamma(\xi_0-\mu u_0)$.
\end{proof}

\begin{remark}
Since we have an existence and uniqueness result in the extended
state space, Theorem~\ref{profCFR} provides an existence and uniqueness
result for \eqref{ABS}, according to Definition~\ref{EUNEW},
whenever we restrict to initial states with $F_0\in\Gamma\PP$.
\end{remark}

However, there are situations where the equality $\SS=\Gamma\PP$
holds, as in the case of the exponential kernel.

\begin{example}
\label{EEXP}
For $a>0$ and $\kappa>0$, consider the kernel
$$\mu(s)=a\e^{-\kappa s},\quad \Omega=\R^+.$$
Since 
$$\mu(t+s)=\e^{-\kappa t}\mu(s),$$
it is apparent that
$$\SS=\SS_0=\big\{F(t)=\e^{-\kappa t}u\,\text{ with }\,u\in V\big\}.$$
In turn,
$$\PP=\big\{\xi(\tau)=\e^{-\kappa \tau}u\,\text{ with }\,u\in V\big\}.$$
Clearly, $\SS=\Gamma\PP$.
\end{example}

\begin{remark}
Incidentally, the above example sheds light on another important issue:
there exist states which are not proper states; in other words,
the inclusion $\PP\subset\V$ is strict (and not even dense).
\end{remark}

In summary, there might be state functions
of the original approach that have no corresponding (proper)
states. Conversely, only the proper states describe
the original problem. In this respect, the state approach is
a more general model, which is able
to describe within the formalism of semigroups also a certain class
of Volterra equations with nonautonomous forcing terms.

Nonetheless, if we start from a proper state, it is reasonable
to expect that the evolution remains confined in the space
of proper states.
To this end, let us define the {\it extended proper state space}
as
$$\H_p=V\times H\times\PP,
$$
which is a normed subspace of $\H$.

\begin{proposition}
If $z\in\H_p$, it follows that $S(t)z\in\H_p$.
\end{proposition}

\begin{proof}
Let $z=(u_0,v_0,\xi_0)\in\H_p$.
Then, $\xi_0=DF$ for some $F\in\SS$. In turn, $F=\Lambda\phi$
for some $\phi\in\AA$.
Denoting as usual $S(t)z=(u(t),v(t),\xi^t)$, and setting
$$\psi^t(s)=u(t-s)\chi_{(0,t)}(s)+u_0\chi_{(t,\ell)}(s)-u(t)$$
and
$$
\phi^t(s)=\phi(s-t)\chi_{(t,\ell)}(s),
$$
the representation formula \eqref{REP} can be equivalently written as
$$\xi^t(\tau)=\int_0^\ell\mu'(\tau+s)\big\{\psi^t(s)+\phi^t(s)\big\}\d s.
$$
By Lemma~\ref{LEMMAPROPER}, in order to prove that $\xi^t\in\PP$,
we are left to show that
$\psi^t+\phi^t\in\AA$.
Indeed, since $\|S(t)z\|_\H\leq\|z\|_\H$,
$$\int_0^\ell\mu(s)\|\psi^t(s)\|_V\d s\leq3M(0)\|z\|_\H.
$$
Therefore, $\psi^t\in L^1_\mu(\Omega;V)\subset\AA$.
Concerning $\phi^t$, we have
$$\int_0^\ell \mu(s)\phi^t(s)\d s=\Lambda\phi(t)\in V,
$$
which yields $\phi^t\in L^1_\mu(\Omega;V)\subset\AA$.
\end{proof}

In particular, from Theorem~\ref{SEMIGROUP}
and Theorem~\ref{EXPthm}, we have the following corollary.

\begin{corollary}
The restriction
$$S_p(t)=S(t)_{|\H_p}:\H_p\to\H_p$$
is a contraction semigroup on $\H_p$.
Assuming also condition \eqref{SUF},
the semigroup $S_p(t)$ is exponentially stable.
\end{corollary}

One might ask whether $\PP$ (and in turn $\H_p$) is a Banach space.
This is true, for instance, for the exponential kernel of Example~\ref{EEXP}.
However, in general, the answer is negative.

\begin{example}
Take the kernel
$$\mu(s)=1-s,\quad\Omega=(0,1).
$$
Let $\CC:[0,1]\to[0,1]$ denote the famous Vitali-Cantor-Lebesgue singular function,
and let $\CC_n$ be the usual approximating sequence of absolutely continuous functions
(cf.\ \cite{HS}).
Consider the sequence
$$F_n(t)=\Big[\int_0^{1-t}\CC_n(s)\d s\Big]\chi_{[0,1]}(t)\,u,
$$
where $u\in V$ is any nonzero vector.
Then, $F_n=\Lambda\phi_n$ with
$$\phi_n(s)=\CC'_n(s)\,u.$$
Setting
$$\xi_n(\tau)=DF_n(\tau)=-\CC_n(1-\tau)\,u,
$$
it is readily seen that $\xi_n\in\V$, and, consequently, $\xi_n\in\PP$.
It also apparent that
$$\lim_{n\to \infty}\xi_n=\xi\quad\text{in }\V,$$
where
$$\xi(\tau)=-\CC(1-\tau)\,u.
$$
However, $\xi\not\in\PP$. Indeed, if not so,
the function
$$F(t)=\Gamma\xi(t)=\Big[\int_0^{1-t}\CC(s)\d s\Big]\chi_{[0,1]}(t)\,u
$$
belongs to $\SS$. Hence, there is
$\phi\in\AA$ such that
$$\Lambda\phi(t)=\int_0^{1-t}\Big(\int_0^s\phi(\sigma)\d \sigma\Big)\d s
=\Big[\int_0^{1-t}\CC(s)\d s\Big]u.
$$
But this implies that
$$\phi(s)=\CC'(s)\,u=0$$
for almost every $s\in(0,1)$. Thus, $F=0$ and $\xi=DF=0$, leading
to a contradiction.
\end{example}
%%%%%%%%%%%%%%%%%%%%%%%%%%%%%%%%%%%%%%%%%%%%

%%%%%%%%%%%%%%%%%%%%%%%%%%%%%%%%%%%%%%%%%%%%
\section{State versus History}
\label{S12}

\noindent
We finally turn to the main issue
that motivated this work: the comparison between
the past history and the state approaches.
We begin to show that each element of $\M$ gives rise to
a proper state, defining the linear map
$$\Pi:\M\to\PP$$
as
$$\Pi\eta(\tau)=-\int_0^\ell\mu'(\tau+s)\eta(s)\d s.$$

\begin{lemma}
\label{lemmaSH}
Let $\eta\in\M$. Then, the vector
$$\Pi\eta(\tau)=-\int_0^\ell\mu'(\tau+s)\eta(s)\d s$$
belongs to $\PP$. Moreover,
$$\|\Pi\eta\|_{\V}\leq\|\eta\|_\M.$$
\end{lemma}

\begin{proof}
Let $\eta\in\M$. Then,
\begin{align*}
\|\Pi\eta(\tau)\|_V^2
&\leq \Big(\int_0^\ell -\mu'(\tau+s)\|\eta(s)\|_V\d s\Big)^2\\
&\leq \int_0^\ell -\mu'(\tau+s)\d s
\int_0^\ell -\mu'(\tau+s)\|\eta(s)\|_V^2\d s\\
&=\mu(\tau)\int_0^\ell -\mu'(\tau+s)\|\eta(s)\|_V^2\d s.
\end{align*}
Therefore,
\begin{align*}
\|\Pi\eta\|_{\V}^2
&\leq\int_0^\ell\d\tau \int_0^\ell -\mu'(\tau+s)\|\eta(s)\|_V^2\d s\\
&=\int_0^\ell  \Big(\int_0^\ell -\mu'(\tau+s)\d\tau\Big)\|\eta(s)\|_V^2\d s\\
&=\|\eta\|^2_{\M}.
\end{align*}
Thus $\Pi\eta\in\V$, and the norm inequality stated above holds
(in fact, equality for $\eta(s)=u$, with $u\in V$).
Since $\M\subset L^1_\mu(\Omega;V)$,
because of the straightforward estimate
$$\int_0^\ell\mu(s)\|\eta(s)\|_V\d s
\leq \sqrt{M(0)}\,\|\eta\|_\M,$$
and
$L^1_\mu(\Omega;V)\subset\AA$, it follows from
Lemma~\ref{LEMMAPROPER} that $\Pi\eta$ is a proper state.
\end{proof}

Rephrasing the lemma, $\Pi\in L(\M,\V)$; namely, $\Pi$ is a
bounded linear operator from $\M$ into $\V$. Moreover
$$\|\Pi\|_{L(\M,\V)}=1.$$

We now clarify the correspondence between $\eta\in\M$ and 
its related proper state $\Pi\eta$.
Letting
$$\bar z=(u_0,v_0,\eta_0)\in\W,\quad
z=(u_0,v_0,\Pi\eta_0)\in\H_p,$$
and denoting
$$\Sigma(t)\bar z=(\bar u(t),\pt \bar u(t),\bar\eta^t),\quad
S_p(t)z=(u(t),\pt u(t),\xi^t),$$
we have the following result.

\begin{proposition}
\label{CORRetaPS}
The equalities
$$u(t)=\bar u(t)\quad\text{and}\quad\xi^t=
\Pi \bar\eta^t
$$
hold for every $t\geq 0$.
\end{proposition}

\begin{proof}
Introduce the function (cf.\ \eqref{REPETA})
$$
\eta^t(s)=
\begin{cases}
u(t)-u(t-s) & 0<s\le t,\\
\eta_0(s-t)+u(t)-u_0 & s>t,
\end{cases}
$$
which solves the Cauchy problem in $\M$
\begin{equation}
\label{ALFA}
\begin{cases}
\displaystyle \ddt \eta^t=T \eta^t+\pt u(t),\\
\noalign{\vskip1mm}
\eta^0=\eta_0.
\end{cases}
\end{equation}
The representation formula \eqref{REP} for $\xi^t$
furnishes
$$
\xi^t(\tau)=\Pi\eta_0(t+\tau)+\mu(\tau)u(t)-\mu(t+\tau)u_0+
\int_0^t\mu'(\tau+s) u(t-s)\d s
=\Pi\eta^t(\tau).
$$
Thus, exploiting \eqref{TEMPO},
$$
\int_0^\ell \xi^t(\tau)\d\tau=
\int_0^\ell \Pi\eta^t(\tau)\d\tau=
\Lambda\eta^t(0)=
\int_0^\ell \mu(s)\eta^t(s)\d s,$$
and, consequently,
\begin{equation}
\label{BETA}
\ptt u+A\Big[u+\int_0^\ell\mu(s)\eta^t(s)\d s\Big]
=\ptt u+A\Big[u+\int_0^\ell \xi^t(\tau)\d\tau\Big]=0.
\end{equation}
Since 
$$u(0)=u_0\quad\text{and}\quad\pt u(0)=v_0,$$
collecting
\eqref{ALFA}-\eqref{BETA} we conclude that
$$(u(t),\pt u(t),\eta^t)=\Sigma(t)\bar z=(\bar u(t),\pt \bar u(t),\bar\eta^t).$$
This finishes the proof.
\end{proof}

Nonetheless, in general, the map $\Pi:\M\to\PP$ is not injective.
This means that two
{\it different} initial histories may entail the {\it same} initial
proper state,
so leading to the same dynamics in the future.

\begin{example}
\label{EXNN}
Let $N\in\N$. Given $a_n>0$ and
$\kappa_N>\ldots>\kappa_1>0$,
consider the kernel
$$\mu(s)=\sum_{n=1}^Na_n\e^{-\kappa_n s},\quad\Omega=\R^+.$$
For $x_m\in\R$ to be determined later, define
$$\eta_0(s)=u,\quad\eta_N(s)=\Big[\sum_{m=1}^Nx_m s^m\Big]u,$$
where $u\in V$ is a fixed nonzero vector.
Clearly, $\eta_0,\eta_N\in \M$ and $\eta_0\neq\eta_N$.
Besides,
$$\Pi\eta_0(\tau)
=\Big[\sum_{n=1}^N a_n \e^{-\kappa_n\tau}\Big]u$$
and
$$\Pi\eta_N(\tau)
=\Big[\sum_{n=1}^N a_nJ_n\e^{-\kappa_n\tau}\Big]u,$$
having set
$$J_n=\kappa_n\sum_{m=1}^N x_m \int_0^\infty s^m \e^{-\kappa_n s} \d s
=\sum_{m=1}^N b_{nm}x_m,$$
with
$$b_{nm}=\frac{m!}{\kappa_n^m}.$$
The determinant of the matrix ${\mathbb B}=\{b_{nm}\}$
is given by
$${\rm det}(\mathbb B)=\prod_{1\leq n\leq N} \frac{n!}{\kappa_n}
\prod_{1\leq m<n\leq N}\Big(\frac1\kappa_n-\frac1\kappa_m\Big)
\neq 0.
$$
Therefore ${\mathbb B}$
is nonsingular, and 
we can choose $\boldsymbol{x}=[x_1,\ldots,x_N]^\top$ to be the (unique) solution
to the linear system
$${\mathbb B}\boldsymbol{x}=[1,\ldots,1]^\top.$$
In which case, $J_n=1$ for all $n$, so that
the equality $\Pi\eta_0=\Pi\eta_N$ holds true.
\end{example}

However, for the kernel of Example~\ref{EXNN}, one can verify
that $\Pi$ maps $\M$ {\it onto} $\PP$. Thus, every proper state is realized 
by a history from $\M$.
On the contrary, the
next example describes a situation where the map $\Pi$ is injective
on $\M$, but $\Pi\M$ is strictly contained in $\PP$, meaning that
all different histories in $\M$ lead to different proper states,
but there are proper states which do not come from histories.
We need first a definition and some preliminary results.

\begin{definition}
A positive sequence $\{\kappa_n\}$, $n\in\N$, is called a {\it M\"untz sequence} if
$\kappa_n\uparrow \infty$ and
$$
\sum_{n=1}^{\infty}\frac 1 {\kappa_n}=\infty.
$$
\end{definition}

Given a function $g\in L^1_{\rm loc}([0,\infty))$
such that $s\mapsto \e^{-\lambda s}g(s)\in L^1(\R^+)$, for some $\lambda>0$,
we denote its (real) Laplace transform
by
$$
\LL g(x)=\int_0^\infty\e^{-x s} g(s)\d s.
$$
A celebrated result due to C.\ M\"untz says that if
$\{\kappa_n\}$ is a M\"untz sequence belonging to the domain of $\LL g$ and
$$
\LL g(\kappa_n)=0,\quad \forall n\in\N,
$$
then $g$ is identically
zero (see \cite{WID}).

The following lemma is standard. A three-line proof is included
for the reader's convenience.

\begin{lemma}
\label{lemmaEXPZERO}
Let $\kappa_n>0$ be strictly increasing, and let $\beta_n\in\R$ be the general term
of an absolutely convergent series.
Consider the function $h:[0,\infty)\to\R$ defined as
$$
h(t)=\sum_{n=1}^\infty \beta_n\e^{-\kappa_n t}.
$$
Then, $h$ is identically zero if and only if $\beta_n=0$ for every $n$.
\end{lemma}

\begin{proof}
One implication is trivial. If $h\equiv 0$, we have the equality
$$0=\int_0^t \e^{\kappa_1 \tau}h(\tau)\d \tau
=\beta_1 t+\sum_{n=2}^\infty \frac{\beta_n}{\kappa_n-\kappa_1}
\Big(1-\e^{-(\kappa_n-\kappa_1) t}\Big),\quad\forall t\geq 0.
$$
The uniform boundedness of the series forces
$\beta_1=0$. Iterate the argument for all $n$.
\end{proof}

We are now ready to provide the aforesaid example.

\begin{example}
Consider the kernel
$$\mu(s)=\sum_{n=1}^\infty a_n\e^{-\kappa_n s},\quad\Omega=\R^+,$$
with $\kappa_n>0$ strictly increasing and $a_n>0$ such that
$$\sum_{n=1}^\infty a_n<\infty.$$
Such a $\mu$ is summable on $\R^+$.
We first observe that if $g\in L^1_\mu(\R^+)$, then
$g\in L^1_{\rm loc}([0,\infty))$ and
$\{\kappa_n\}$ belongs to the domain of $\LL g$.
Let us extend in the obvious way the map $\Pi$ to the domain
$$\M_\star=\Big\{\eta\in L^1_\mu(\Omega;V)\,:\,
\tau\mapsto -\int_0^\ell\mu'(\tau+s)\eta(s)\d s\in \V\Big\}
$$
(we keep calling $\Pi$ such an extension).
Note that $\M\subset\M_\star\subset\AA$, and from Lemma~\ref{LEMMAPROPER}
we learn that $\Pi\M_\star\subset\PP$.
Given $\eta\in\M_\star$ and $w\in V^*$,
we consider the duality product
$$g_w(s)=\l\eta(s),w\r\in L^1_\mu(\R^+).$$
In view of \eqref{TEMPO},
$$\Pi\eta=0\quad\Leftrightarrow\quad \Lambda\eta=0.
$$
But $\Lambda\eta=0$ if and only if
$$
\sum_{n=1}^\infty \beta_n(w)\e^{-\kappa_n t}=0,\quad\forall t\geq 0,\,\forall w\in V^*,
$$
having set $\beta_n(w)= a_n\,\LL g_w(\kappa_n)$.
Moreover,
$$
\sum_{n=1}^\infty |\beta_n(w)|
\leq\sum_{n=1}^\infty a_n\int_0^\infty\e^{-\kappa_n s} |g_w(s)|\d s
=\|g_w\|_{L^1_\mu(\R^+)}<\infty.
$$
Hence, from Lemma~\ref{lemmaEXPZERO},
the above equality is true if and only if
$$
\LL g_w(\kappa_n)=0,\quad \forall n\in\N,\,\forall w\in V^*.
$$
Therefore, if $\{\kappa_n\}$ is a M\"untz sequence,
$$\Pi\eta=0\quad\Leftrightarrow\quad g_w= 0,\,\,\forall w\in V^*
\quad\Leftrightarrow \quad \eta=0.
$$
In which case, the map $\Pi$ is injective on $\M_\star$.
Accordingly, to conclude that $\Pi\M$ is a proper subset of $\PP$
we have to show that the inclusion $\M\subset\M_\star$ is strict.
This is obtained, for instance, by looking at the elements
$$\eta(s)=\e^{\sigma\kappa_1 s}u,$$
where $\sigma\in[\frac12,1)$ and $u\in V$ is any nonzero vector.
The details are left to the reader.
\end{example}
%%%%%%%%%%%%%%%%%%%%%%%%%%%%%%%%%%%%%%%%%%%%

%%%%%%%%%%%%%%%%%%%%%%%%%%%%%%%%%%%%%%%%%%%%

%%%%%%%%%%%%%%%%%%%%%%%%%%%%%%%%%%%%%%%%%%%%


\begin{thebibliography}{99}

\bibitem{BOL1}
{\au L.~Boltzmann},
{\ti Zur Theorie der elastischen Nachwirkung},
{\jou Wien.\ Ber.}
\no{70}{275--306}{1874}

\bibitem{BOL2}
{\au L.~Boltzmann},
{\ti Zur Theorie der elastischen Nachwirkung},
{\jou Wied.\ Ann.}
\no{5}{430--432}{1878}

\bibitem{BO}
{\au S.~Breuer, E.T.~Onat},
{\ti On recoverable work in linear viscoelasticity},
{\jou Z.\ Angew.\ Math.\ Phys.}
\no{15}{13--21}{1964}

\bibitem{CP}
{\au V.V.~Chepyzhov, V.~Pata},
{\ti Some remarks on stability of semigroups arising from
linear viscoelasticity}
{\jou Asymptot.\ Anal.}
\no{50}{269--291}{2006}

\bibitem{COL}
{\au B.D.~Coleman},
{\ti Thermodynamics of materials with memory},
{\jou Arch.\ Ration.\ Mech.\ Anal.}
\no{17}{1--45}{1964}

\bibitem{CM1}
{\au B.D.~Coleman, V.J.~Mizel},
{\ti Norms and semi-groups in the theory of fading memory},
{\jou Arch.\ Ration.\ Mech.\ Anal.}
\no{23}{87--123}{1967}

\bibitem{CM2}
{\au B.D.~Coleman, V.J.~Mizel},
{\ti On the general theory of fading memory},
{\jou Arch.\ Ration.\ Mech.\ Anal.}
\no{29}{18--31}{1968}

\bibitem{CN}
{\au B.D.~Coleman, W.~Noll},
{\ti Foundations of linear viscoelasticity},
{\jou Rev.\ Modern Phys.}
\no{33}{239--249}{1961}

\bibitem{DAF1}
{\au C.M.~Dafermos},
{\ti Asymptotic stability in viscoelasticity},
{\jou Arch.\ Ration.\ Mech.\ Anal.}
\no{37}{297--308}{1970}

\bibitem{DAY1}
{\au W.A.~Day},
{\ti Reversibility, recoverable work and free energy in linear  
viscoelasticity},
{\jou Quart.\ J.\ Mech.\ Appl.\ Math.}
\no{23}{1--15}{1970}

\bibitem{DAY2}
{\au W.A.~Day},
{\bk The thermodynamics of simple materials with fading memory},
\eds{Springer}{New York}{1972}

\bibitem{DPD1}
{\au G.~Del Piero, L.~Deseri},
{\ti Monotonic, completely monotonic and exponential relaxation  
functions in linear viscoelasticity},
{\jou Quart.\ Appl.\ Math.}
\no{53}{273--300}{1995}

\bibitem{DPD2}
{\au G.~Del Piero, L.~Deseri},
{\ti On the concepts of state and free energy in linear  
viscoelasticity},
{\jou Arch.\ Ration.\ Mech.\ Anal.}
\no{138}{1--35}{1997}

\bibitem{DFG}
{\au L.~Deseri, M.~Fabrizio, M.J.~Golden},
{\ti The concept of minimal state in viscoelasticity: new free energies  
an applications to PDEs},
{\jou Arch.\ Ration.\ Mech.\ Anal.}
\no{181}{43--96}{2006}

\bibitem{DGG}
{\au L.~Deseri, G.~Gentili, M.J.~Golden},
{\ti An explicit formula for the minimum free energy in linear  
viscoelasticity},
{\jou J.\ Elasticity}
\no{54}{141--185}{1999}

\bibitem{FGM0}
{\au M.~Fabrizio, C.~Giorgi, A.~Morro},
{\ti Minimum principles, convexity,
and thermodynamics in linear viscoelasticity},
{\jou Continuum Mech.\ Thermodyn.}
\no{1}{197--211}{1989}

\bibitem{FGM1}
{\au M.~Fabrizio, C.~Giorgi, A.~Morro},
{\ti Free energies and dissipation properties for systems with memory},
{\jou Arch.\ Ration.\ Mech.\ Anal.}
\no{125}{341--373}{1994}

\bibitem{FGM2}
{\au M.~Fabrizio, C.~Giorgi, A.~Morro},
{\ti Internal dissipation, relaxation property and free energy in  
materials with fading memory},
{\jou J.\ Elasticity}
\no{40}{107--122}{1995}

\bibitem{FG1}
{\au M.~Fabrizio, M.J.~Golden},
{\ti Maximum and minimum free energies for a linear viscoelastic  
material},
{\jou Quart.\ Appl.\ Math.}
\no{60}{341--381}{2002}

\bibitem{FG2}
{\au M.~Fabrizio, M.J.~Golden},
{\ti Minimum free energies for materials with finite memory},
{\jou J.\ Elasticity}
\no{72}{121--143}{2003}

\bibitem{FL}
{\au M.~Fabrizio, B.~Lazzari},
{\ti On the existence and asymptotic stability of
solutions for linear viscoelastic solids},
{\jou Arch.\ Ration.\ Mech.\ Anal.}
\no{116}{139--152}{1991}

\bibitem{FL1}
{\au M.~Fabrizio, B.~Lazzari},
{\ti Stability and free energies in linear viscoelasticity},
{\jou Matematiche (Catania)}
\no{62}{175--198}{2007}

\bibitem{FM0}
{\au M.~Fabrizio, A.~Morro},
{\ti Viscoelastic relaxation functions compatible with thermodynamics},
{\jou J.\ Elasticity}
\no{19}{63--75}{1988}

\bibitem{FM}
{\au M.~Fabrizio, A.~Morro},
{\bk Mathematical problems in linear viscoelasticity},
\eds{SIAM Studies in Applied Mathematics no.12,  
SIAM}{Philadelphia}{1992}

\bibitem{FIC1}
{\au G.~Fichera},
{\ti Analytic problems of hereditary phenomena in materials with  
memory},
in ``Corso CIME''
(Bressanone, 1977),
\eds{pp.111-169, Liguori}{Napoli}{1979}

\bibitem{FIC2}
{\au G.~Fichera},
{\ti Avere una memoria tenace crea gravi problemi},
{\jou Arch.\ Ration.\ Mech.\ Anal.}
\no{70}{101--112}{1979}

\bibitem{GEN}
{\au G.~Gentili},
{\ti Maximum recoverable work, minimum free energy and state space in  
linear viscoelasticity},
{\jou Quart.\ Appl.\ Math.}
\no{60}{152--182}{2002}

\bibitem{GPR}
{\au C.~Giorgi, J.E.~Mu\~noz Rivera, V.~Pata},
{\ti Global attractors for a semilinear hyperbolic equation in  
viscoelasticity},
{\jou J.\ Math.\ Anal.\ Appl.}
\no{260}{83--99}{2001}

\bibitem{GG}
{\au J.M.~Golden, G.A.C.~Graham},
{\bk Boundary value problems in linear viscoelasticity},
\eds{Springer}{New York}{1988}

\bibitem{GRA1}
{\au D.~Graffi},
{\ti Sui problemi della eredit\`a lineare},
{\jou Nuovo Cimento}
\no{5}{53--71}{1928}

\bibitem{GRA2}
{\au D.~Graffi},
{\it Sopra alcuni fenomeni ereditari dell'elettrologia},
{\jou Rend.\ Istit.\ Lombardo Sc.\ Lett.}
\no{68-69}{124--139}{1936}

\bibitem{GRA3}
{\au D.~Graffi},
{\ti Sull'espressione analitica di alcune grandezze termodinamiche nei  
materiali con memoria},
{\jou Rend.\ Sem.\ Mat.\ Univ.\ Padova}
\no {68}{17--29}{1982}

\bibitem{GRA4}
{\au D.~Graffi},
{\ti On the fading memory},
{\jou Appl.\ Anal.}
\no {15}{295--311}{1983}

\bibitem{GF}
{\au D.~Graffi, M.~Fabrizio},
{\ti Sulla nozione di stato per materiali viscoelastici di tipo  
``rate"},
{\jou Atti Accad.\ Lincei Rend.\ Fis.}
\no{83}{201--208}{1989}

\bibitem{GP-Terreni}
{\au M.~Grasselli, V.~Pata},
{\ti Uniform attractors of nonautonomous systems with memory},
in ``Evolution Equations, Semigroups and Functional Analysis''
(A.~Lorenzi and B.~Ruf, Eds.),
\eds{pp.155--178, Progr.\ Nonlinear Differential Equations
Appl.\ no.50, Birkh\"{a}user}{Boston}{2002}

\bibitem{GRIV}
{\au A.E.~Green, R.S.~Rivlin},
{\ti The mechanics of nonlinear materials with memory},
{\jou Arch.\ Ration.\ Mech.\ Anal.}
\no{1}{1--21}{1957-58}

\bibitem{GS}
{\au M.E.~Gurtin, E.~Sternberg},
{\ti On the linear theory of viscoelasticity},
{\jou Arch.\ Ration.\ Mech.\ Anal.}
\no{11}{291--356}{1962}

\bibitem{HS}
{\au E.~Hewitt, K.~Stromberg},
{\bk Real and abstract analysis},
\eds{Springer-Verlag}
{New York}{1965}

\bibitem{II}
{\au M.G.~Ianniello, G.~Israel},
{\ti Boltzmann's concept of ``Nachwirkung" and the ``mechanics of  
heredity"},
in ``Proceedings of the International
Symposium on Ludwig Boltzmann"
(G.~Battimelli, M.G.~Ianniello and O.~Kresten, Eds.)
\eds{pp.113-133, Verlag der {\"O}sterreichischen
Akademie der Wissenschaften}{Wien}{1993}

\bibitem{KM}
{\au H.~K\"onig, J.~Meixner},
{\ti Lineare Systeme und lineare Transformationen},
{\jou Math.\ Nachr.}
\no{19}{265--322}{1958}

\bibitem{LF}
{\au M.J.~Leitman, G.M.C.~Fisher},
{\ti The linear theory of viscoelasticity},
in ``Handbuch der Physik" vol.VIa/3
(S.~Fl\"ugge, Ed.),
\eds{pp.1--123, Springer}{Berlin-Heidelberg-New York}{1973}

\bibitem{LZ0}
{\au Z.~Liu, S.~Zheng},
{\ti On the exponential stability
of linear viscoelasticity and thermoviscoelasticity},
{\jou Quart.\ Appl.\ Math.}
\no{54}{21--31}{1996}

\bibitem{RIV}
{\au J.E.~Mu\~noz Rivera},
{\ti Asymptotic behaviour in linear viscoelasticity},
{\jou Quart.\ Appl.\ Math.}
\no{52}{629--648}{1994}

\bibitem{NOL}
{\au W.~Noll},
{\ti A new mathematical theory of simple materials},
{\jou Arch.\ Ration.\ Mech.\ Anal.}
\no{48}{1--50}{1972}

\bibitem{PAT}
{\au V.~Pata},
{\ti Exponential stability in linear viscoelasticity},
{\jou Quart.\ Appl.\ Math.}
\no{64}{499--513}{2006}

\bibitem{PAZ}
{\au A.~Pazy},
{\bk Semigroups of linear operators and applications to partial  
differential
equations},
\eds{Springer-Verlag}{New York}{1983}

\bibitem{RHN}
{\au M.~Renardy, W.J.~Hrusa, J.A.~Nohel},
{\bk Mathematical problems in viscoelasticity},
\eds{Harlow John Wiley \& Sons}{New York}{1987}

\bibitem{VOL1}
{\au V.~Volterra},
{\ti Sur les \'equations int\'egro-diff\'erentielles et leurs  
applications},
{\jou Acta Math.}
\no{35}{295--356}{1912}

\bibitem{VOL2}
{\au V.~Volterra},
{\bk Le\c{c}ons sur les fonctions de lignes},
\eds{Gauthier-Villars}{Paris}{1913}

\bibitem{WID}
{\au D.V.~Widder},
{\bk The Laplace transform},
\eds{Princeton University Press}{Princeton}{1941}

\end{thebibliography}
\end{document}